\documentclass[12pt,english]{article} 
\usepackage{mathptmx}
\usepackage{geometry}
\usepackage{color}
\usepackage{xcolor}
\usepackage{array}
\usepackage{bbding}
\usepackage{multirow}
\usepackage[fleqn]{amsmath}
\usepackage{amssymb}
\usepackage{amsthm}
\usepackage{graphicx}
\usepackage{natbib}
\usepackage{lscape}
\usepackage[hyphens]{url}
\usepackage[hidelinks]{hyperref}
\usepackage{comment}
\usepackage[title]{appendix}
\usepackage{longtable}
\usepackage{mathtools}
\usepackage{chngcntr}
\usepackage{authblk}
\usepackage{babel}
\usepackage{dsfont}
\usepackage{subcaption}
\usepackage{fancyhdr}
\usepackage{cleveref}
\usepackage{abstract}

\makeatletter
\allowdisplaybreaks
\usepackage[utf8]{inputenc}
\usepackage[T1]{fontenc}
\usepackage{newpxtext, newpxmath}

\date{}

\fancypagestyle{plain}{%
	\fancyhf{} 
	
	\chead{\footnotesize \textcolor{gray}{\textit{\shortauthors -- \runningtitle}}}
	\rhead{\footnotesize \textcolor{gray}{\thepage}}
}
\makeatletter
\def\blfootnote{\xdef\@thefnmark{}\@footnotetext}
\makeatother

\makeatletter
\def\titlepageext{
	\begin{center}	
	{\parindent0pt
		\rule{0.9\textwidth}{1pt}
		\begin{minipage}[t]{0.25\textwidth}
			\small {\it Keywords:}\\
			\keyword
		\end{minipage}%
		\hspace{3mm}
		\begin{minipage}[t]{0.6\textwidth}
			\small \abstract
		\end{minipage}%
		
		\rule{0.9\textwidth}{2pt}
	}
	\end{center}

	\blfootnote{* Corresponding author. E-mail address: \href{mailto:\corresemail}{\corresemail}.}
}
\makeatother

\usepackage[mathlines, pagewise]{lineno}
\usepackage{etoolbox} 

\newcommand*\linenomathpatchAMS[1]{%
	\expandafter\pretocmd\csname #1\endcsname {\linenomathAMS}{}{}%
	\expandafter\pretocmd\csname #1*\endcsname{\linenomathAMS}{}{}%
	\expandafter\apptocmd\csname end#1\endcsname {\endlinenomath}{}{}%
	\expandafter\apptocmd\csname end#1*\endcsname{\endlinenomath}{}{}%
}

\expandafter\ifx\linenomath\linenomathWithnumbers
\let\linenomathAMS\linenomathWithnumbers
\patchcmd\linenomathAMS{\advance\postdisplaypenalty\linenopenalty}{}{}{}
\else
\let\linenomathAMS\linenomathNonumbers
\fi

\linenomathpatchAMS{gather}
\linenomathpatchAMS{multline}
\linenomathpatchAMS{align}
\linenomathpatchAMS{alignat}
\linenomathpatchAMS{flalign}

\nolinenumbers

\newtheorem{theorem}{Theorem}

\usepackage{mathtools}

\usepackage{longtable}
\usepackage{tabularx}

\usepackage{tikz}
\usepackage{pgfplots}
\usetikzlibrary{shapes.geometric}
\usetikzlibrary{arrows.meta}
\usetikzlibrary{overlay-beamer-styles}
\usetikzlibrary{calc,3d}
\usetikzlibrary{decorations.pathreplacing}
\usetikzlibrary{math}
 \usetikzlibrary{patterns}

\usepackage{graphicx}
\usepackage{multirow}
\usepackage{hhline}
\usepackage{wasysym}
\def\Halmos{\mbox{\quad$\square$}}
\makeatletter
\newcommand{\UE}[1]{#1^{\mathrm{UE}}}
\newcommand{\SO}[1]{#1^{\mathrm{SO}}}
\newcommand{\NS}[1]{#1^{\mathrm{NS}}}
\newcommand{\TLC}[1]{#1^{\mathrm{TLC}}}
\newcommand{\SWH}[1]{#1^{\mathrm{SWH}}}
\newcommand{\CS}[1]{#1^{\mathrm{CS}}}
\makeatother

\usepackage{amsthm}
\theoremstyle{definition}
\theoremstyle{definition}

\newtheorem{definition}{Definition}[section]
\newtheorem{lemma}{Lemma}[section]
\newtheorem{assumption}{Assumption}[section]
\newtheorem{proposition}{Proposition}[section]

\crefname{section}{Section}{Sections}
\crefname{subsection}{Section}{Sections}

\crefname{equation}{Eq.}{Eqs.}
\crefname{figure}{Figure}{Figures}
\crefname{subfigure}{Figure}{Figures}
\crefname{table}{Table}{Tables}

\crefname{thm}{Theorem}{Theorems}
\crefname{cor}{Corollary}{Corollary}
\crefname{dfn}{Definition}{Definition}
\crefname{lem}{Lemma}{Lemmas}
\crefname{asm}{Assumption}{Assumption}
\crefname{pro}{Proposition}{Proposition}
\crefname{cnj}{Conjecture}{Conjecture}

\Crefname{section}{Section}{Sections}
\Crefname{subsection}{Section}{Sections}
\Crefname{equation}{Eq.}{Eqs.}
\Crefname{figure}{Figure}{Figures}
\Crefname{subfigure}{Figure}{Figures}
\Crefname{table}{Table}{Tables}
\Crefname{thm}{Theorem}{Theorems}
\Crefname{cor}{Corollary}{Corollary}
\Crefname{dfn}{Definition}{Definition}
\Crefname{lem}{Lemma}{Lemmas}
\Crefname{asm}{Assumption}{Assumption}
\Crefname{pro}{Proposition}{Proposition}
\Crefname{cnj}{Conjecture}{Conjecture}

\DeclareMathOperator*{\argmin}{argmin}

\usepackage[ISO]{diffcoeff}
\usepackage{amsmath,amssymb,bm}
\newcommand{\Vt}[1]{\bm{#1}}

\newcommand{\Vtf}{\bm{f}}

\newcommand{\Vtp}{\bm{p}}
\newcommand{\Vtq}{\bm{q}}
\newcommand{\Vtr}{\bm{r}}

\newcommand{\Vtw}{\bm{w}}

\newcommand{\VtA}{\bm{A}}

\newcommand{\VtP}{\bm{P}}
\newcommand{\VtQ}{\bm{Q}}

\newcommand{\Vtlambda}{\bm{\lambda}}
\newcommand{\Vtmu}{\bm{\mu}}

\newcommand{\ClH}{\mathcal{H}}
\newcommand{\ClI}{\mathcal{I}}

\newcommand{\ClK}{\mathcal{K}}

\newcommand{\ClN}{\mathcal{N}}

\newcommand{\ClT}{\mathcal{T}}

\title{A Paradox of Telecommuting and Staggered Work Hours in the Bottleneck Model}

\def\shortauthors{Sakai et al.}
\def\runningtitle{Paradox of Telecommuting and Staggered Work Hours}

\author[a,$\ast$]{Takara Sakai}
\author[b]{Takashi Akamatsu}
\author[c]{Koki Satsukawa}
\affil[a]{Department of Civil and Environmental Engineering, Tokyo Institute of Technology, Tokyo, 152-8552, Japan}
\affil[b]{Graduate School of Information Sciences, Tohoku University, Miyagi, 980-8579, Japan}
\affil[c]{Institute of Transdisciplinary Sciences for Innovation, Kanazawa University, Ishikawa, 920-1192, Japan.}
\def\corresemail{sakai.t.av@m.titech.ac.jp}

\def\abstract{We study the long- and short-term effects of telecommuting (TLC), staggered work hours (SWH), and their combined scheme on peak-period congestion and location patterns.
In order to enable a unified comparison of the schemes' long- and short-term effects, we develop a novel equilibrium analysis approach that consistently synthesizes the long-term equilibrium (location and percentage of telecommuting choice) and short-term equilibrium (preferred arrival time and departure time choice).
By exploiting their special mathematical structures similar to optimal transport problems, we derive the closed-form solution to the long- and short-term equilibrium while explicitly considering their interaction. 
These closed-form solutions elucidate the discrepancies between the effects of each scheme and uncover a paradoxical finding: the introduction of SWH, in conjunction with TLC, may increase the total commuting costs compared to the scenario with only TLC, without yielding any improvement in worker utility.}

\def\keyword{bottleneck model \\ corridor network \\ telecommuting \\ staggered work hours \\ departure time choice \\ location choice}

\begin{document}
\maketitle
\titlepageext

\section{Introduction}
\subsection{Background and Purpose}
\noindent
Telecommuting (TLC) and staggered work hours (SWH) are flexible forms of work that have been spread over recent decades. These flexible working styles relax temporal and spatial constraints on working, which will affect short-term and long-term phenomena in cities, such as commuting traffic and residential location patterns. Moreover, many companies have recently adopted these working styles to reduce unnecessary contact with others owing to the COVID-19 pandemic. 
Therefore, it is increasingly important for urban and transportation planning to take into account the impact of these changes in working styles on urban economics.
	\par
  Regarding the effects of TLC have been studied extensively since the seminal work of~\citet{Nilles1988-fd}.
As for effects on location patterns, many studies have shown that TLC may drive workers away from the CBD, thus saving office space and changing the workers' location choice farther from the CBD, leading to a more spread-out city~\citep{Janelle1986-iq,Nilles1991-hy,Rhee2009-yp}.
Regarding the effects of TLC on commuting traffic patterns, they have demonstrated that it can effectively reduce the overall demand for commuting, resulting in congestion reduction, energy saving, and air quality improvement~\citep{Safirova2002-ks,Delventhal2022-ix}.
However, the existing studies have primarily modeled commuting traffic using a static framework, which has resulted in a limited understanding of the effects of TLC on within-day congestion dynamics.
\par
In contrast to studies on TLC, some studies focusing on the impacts of SWH have adopted dynamic frameworks to model commuting traffic (e.g., ~\citealp{Ben-Akiva1984-il,Mun2006-vi,Fosgerau2017-ks,Takayama2015-dp}). 
These investigations have shown that SWH helps temporal dispersion of traffic demand, resulting in congestion reduction during peak periods. However, few studies have comprehensively analyzed the long-term effects of SWH on location patterns.
\par
As demonstrated earlier, previous studies provide valuable insights into TLC and SWH; however, there is a lack of knowledge regarding the short-term effect of TLC on within-day commuting patterns and the long-term effect of SWH on location patterns.
Moreover, existing theoretical studies have primarily examined each scheme in isolation, thus providing a limited understanding of the combined implementation of TLC and SWH and their comprehensive effects on both short- and long-term phenomena. 
Considering the growing implementation of TLC and SWH in practice, we need to systematically investigate their effects within the context of urban and transportation planning.
To this end, a comprehensive analytical model is required to capture the short- and long-term consequences of TLC, SWH, and their combined scheme. 
However, to the best of our knowledge, no existing theoretical framework encompasses a comprehensive description of these effects.
\begin{figure}[tbp]
  \center
   \includegraphics[clip, width=0.50\columnwidth]{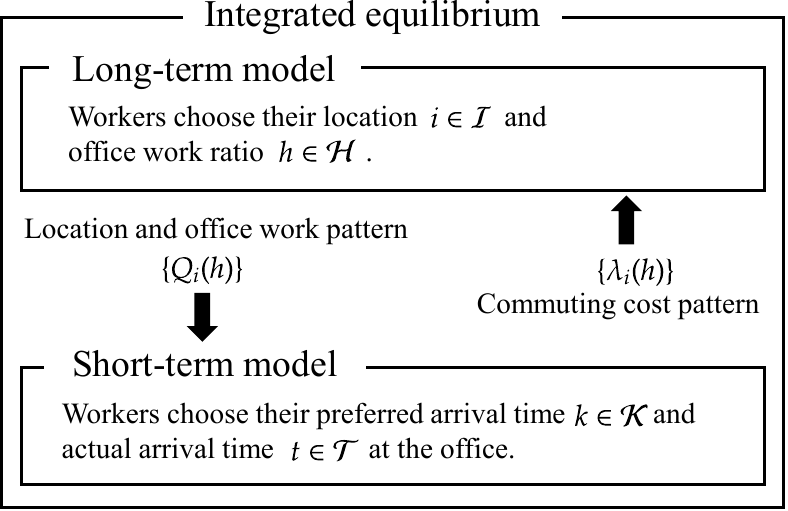}
  \caption{Integrated Equilibrium Model; if $\ClH=\{1\}$ and $\|\ClK\|=1$, then this model describes the equilibrium state without TLC and SWH, respectively.}
  \label{fig:Strucure_IntModel}
\end{figure}
\par
This study aims to clarify the long and short-term effects of TLC, SWH, and their combined scheme on peak-period congestion and location patterns.
To this end, keeping the model as compact as possible, we extend the standard bottleneck model to allow each worker to choose their residential location, office work ratio, and preferred arrival time.
The office work ratio and the preferred arrival time represent the percentage of office working days with commuting per unit term and the official work start time at the office place, respectively, which are introduced for modeling TLC and SWH.
Our model consistently integrates the long- and short-term equilibrium using the interdependence of commuting demand and cost on long- and short-term models.
\Cref{fig:Strucure_IntModel} illustrates the structure of our integrated model (the detailed formulation of the integrated model is described in the later section).
In the long-term model, workers choose their residential location $i \in \ClI$ and office work ratio $h \in \ClH \subseteq [0, 1]$ for given commuting cost patterns $\{ \lambda_{i}(h) \}$ determined by the short-term equilibrium. 
Here, $\ClI$ and $\ClH$ represent the sets of residential locations and office work ratios, respectively. 
In the short-term model, workers choose their preferred arrival time $k \in \ClK$ and actual arrival time $t \in \ClT$ for given commuting demand patterns $\{ Q_{i}(h) \}$ determined by the long-term equilibrium. 
Here, $\ClK$ and $\ClT$ represent the sets of preferred arrival times and actual arrival times, respectively.
Additionally, the integrated equilibrium model encompasses the equilibrium states without TLC or SWH. Specifically, when the set of office work ratios $\ClH$ is $\{1\}$, our model describes the equilibrium state without TLC. 
Similarly, when the size of the set of preferred arrival times $\ClK$ is $1$, our model represents the equilibrium state without SWH. Thus, our framework facilitates a comprehensive comparison of the long- and short-term effects of these schemes.
\par
In order to compare the effects of each scheme, we derive the general closed-form solution to the integrated equilibrium model by a step-by-step approach.
We first derive the short-term equilibrium by the queue replacement principle (QRP) approach~\citep{Akamatsu2021-zg,Fu2022-nl,Sakai2022-ay}. 
The QRP approach is a technique for deriving equilibrium states using social optimal states.
Because we can derive the solution to the short-term optimal problem by utilizing the theory of optimal transport~\citep{Rachev1998-bb}, we obtain the closed-form solution to the short-term equilibrium problem by applying the QRP approach.
This closed-form solution is represented as a function of the demand pattern $\{ Q_{i}(h) \}$ determined by long-term equilibrium.
Using this functional relationship, we then derive the closed-form long-term equilibrium.
Finally, by combining the long- and short-term equilibrium states, we obtain the general closed-form solution for the integrated equilibrium model.
\par
This general solution enables us to derive the equilibrium without TLC or SWH. By comparing these equilibrium states, we demonstrate the positive impact of introducing TLC or SWH on workers' utility, as anticipated. Surprisingly, we find that introducing SWH in addition to TLC does not enhance workers' utility compared to a scenario with TLC alone. Furthermore, we discover that this may lead to higher total commuting costs.

\subsection{Related Literature and Contributions}
The concept of TLC was first proposed by~\citet{Nilles1988-fd} and has been studied to understand its impacts on residential location patterns, traffic patterns, and other social phenomena.
\citet{Higano1990-as}, \citet{Lund1994-rj}, and \citet{Safirova2002-ks} have shown model-based analysis for the long-term impacts of TLC on residential location patterns.
\citet{Delventhal2022-ix} studied how cities change when there were a permanent increase in TLC using a quantitative model.
Although these analyses have included the effects on commuting traffic, most have not considered within-day traffic congestion. 
The impacts of TLC on dynamic traffic congestion have been measured by mainly empirical studies\footnote{Numerous empirical studies that have analyzed the effects of introducing TLC  during the COVID-19 pandemic~\citep{Loo2022-xe,Wohner2022-rg,Tahlyan2022-wj,Zhu2022-yq}.}\citep{Ory2006-oj,Asgari2018-cr,Lachapelle2018-zl,Shabanpour2018-el}.
The limited number of theoretical studies examining the impact of TLC on dynamic congestion include \citet{Gubins2011-qk} and \citet{WADA2023-xm}, who conducted welfare analyses utilizing a single bottleneck model.
\par
The effects of SWH on traffic congestion have been analyzed by many studies since the seminal work of \citet{Henderson1981-yb}.
\citet{Ben-Akiva1984-il,Mun2006-vi,Fosgerau2017-ks}, and \citet{Takayama2015-dp} were the most successful in considering both SWH effects and peak-period traffic congestion.
They showed the basic properties of the equilibrium with SWH and welfare analysis.
The properties of bottleneck congestion with SWH can also be characterized from the research of the bottleneck models that consider the heterogeneity of preferred arrival times (e.g., \citealp{Hendrickson1981-cu,Daganzo1985-ls,Lindsey2004-aw,Lindsey2019-zs}).
Especially first-in-first-work (FIFW) principle shown in these studies helps us to understand the equilibrium flow pattern.
\par
Existing studies have independently investigated the effects of TLC and SWH on commuting traffic patterns, yielding valuable theoretical insights for each scheme. 
However, there is a lack of knowledge regarding the short-term effect of TLC on within-day commuting patterns and the long-term effect of SWH on location patterns.
Moreover, due to the separate focus of existing theoretical studies focused on each scheme, there is limited understanding of how their combined scheme impacts both short- and long-term phenomena.
\par
In contrast, our contributions are summarized as follows:
\begin{itemize}
  \setlength{\itemsep}{0pt}
	\setlength{\parskip}{0pt}      
	\setlength{\itemindent}{0pt}   
	\setlength{\labelsep}{5pt} 
  \item We develop a comprehensive framework for analyzing the integrated equilibrium model, which consistently integrates the long-term equilibrium (location and office work ratio choices) and short-term equilibrium (preferred and actual arrival times choices).
  \item We propose a novel approach based on the theory of optimal transport and the bottleneck decomposition technique to analyze the integrated equilibrium model.
  This approach enables us to derive the closed-form solution to the long- and short-term equilibrium while explicitly considering their interaction.
  \item Through our analysis, we clarify the distinctions between the long- and short-term effects of TLC, SWH, and their combined scheme, and uncover a paradox that arises from the combined scheme.
\end{itemize}

The remainder of this paper is organized as follows: 
Section 2 presents the formulation of the integrated equilibrium model. 
In Section 3, we derive the solution to the short-term problem and explore the theoretical properties of the equilibrium flow pattern. 
Section 4 examines the long-term equilibrium by utilizing the solution to the short-term problem and derives the closed-form solution to the integrated equilibrium model. 
Section 5 demonstrates the effects of TLC and SWH on traffic congestion and residential patterns. 
Finally, in Section 6, we conclude the study and discuss potential avenues for future research.

\begin{figure}[tbp]
  \center
   \includegraphics[clip, width=0.50\columnwidth]{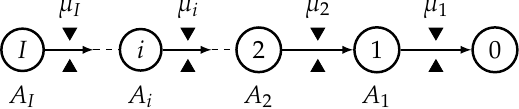}
   \caption{Corridor Network with $I$ Origin and the Single Destination $0$.}
   \label{fig:CorridorNetwork}
\end{figure}

\section{Model}
This section formulates the integrated equilibrium model.
First, we describe the network and worker settings.
Subsequently, we formulate the short- and long-term models, which are part of the integrated equilibrium model.
Finally, by using these models, we formally define the integrated equilibrium model. 
Appendix \ref{sec:List_Notation} presents the list of notation frequently used in this paper.

\subsection{Network and Workers}
Consider a freeway corridor that connects $I$ residential locations to a central business district (CBD), as shown in~\Cref{fig:CorridorNetwork}.
The residential locations are indexed sequentially from the CBD.
We denote the set of locations by $\ClI \equiv \{1,\ldots, I \}$.
Each residential location $i \in \ClI$ is endowed with $A_{i}$ unit of land.
Each link $i \in \ClI \equiv \{1,\ldots, I \}$ connecting from the origin (location) $i$ has a bottleneck with a finite capacity $\mu_{i}$.
The queue evolution and associated queueing delay are modeled by the standard point queue model with the first-in-first-principle.
The free-flow travel time of link $i$ is denoted by $f_{i}$.
We assume $\mu_{i} > \mu_{i+1}$ for all $i \in \ClI \setminus \{ I \}$.
\par
The firms located in the CBD are homogeneous and have the same productivity.\footnote{
  This model does not describe the temporal agglomeration economics of the office work compared to \citet{Henderson1981-yb,Takayama2015-dp}, because we focus on the bottleneck congestion with a single peak-period.
  In single peak bottleneck congestion, the CBD arrival flow rate is constant even if the queueing pattern is changed, which is equal to the bottleneck capacity.
  Therefore, the number of workers in the office at a given instant is not affected by the queueing pattern.
}
The firms allow workers to work remotely and flexibly.
Specifically, the workers can choose the ``office work ratio'' and ``preferred arrival time.''
The office work ratio $h \in \ClH \equiv [0, 1]$ is the  percentage of office working days with commuting per unit term.\footnote{
By imposing the upper and lower bound of the office work ratio $h$ (i.e., $\ClH \equiv [\underline{h}, \overline{h}]$, where $0<\underline{h}< \overline{h}<1$), we can model a situation where firms require workers to come to work part of the time (e.g., 2 or 3 days per week).
In order to simplify the analysis, this study assumes that the upper and lower bounds of the office work ratio are $1$ and $0$, respectively.}
The firms pay the wage $\theta^{O}$ to office workers and the wage $\theta^{R}$ to remote workers (telecommuters) per day, where $\theta^{O}$ and $\theta^{R}$ are given parameters and satisfy $\theta^{O} > \theta^{R}$. 
The firms set $K$ official work start times $\{ t_{1}, ..., t_{K} \}$, indicating the preferred arrival times for the workers.
The difference between preferred arrival times is denoted by $d_{k} = t_{k+1} - t_{k}$ for all $k \in \ClK \setminus \{ K \}$.
\par
\begin{figure}[tbp]
  \center
     \includegraphics[clip, width=0.55\columnwidth]{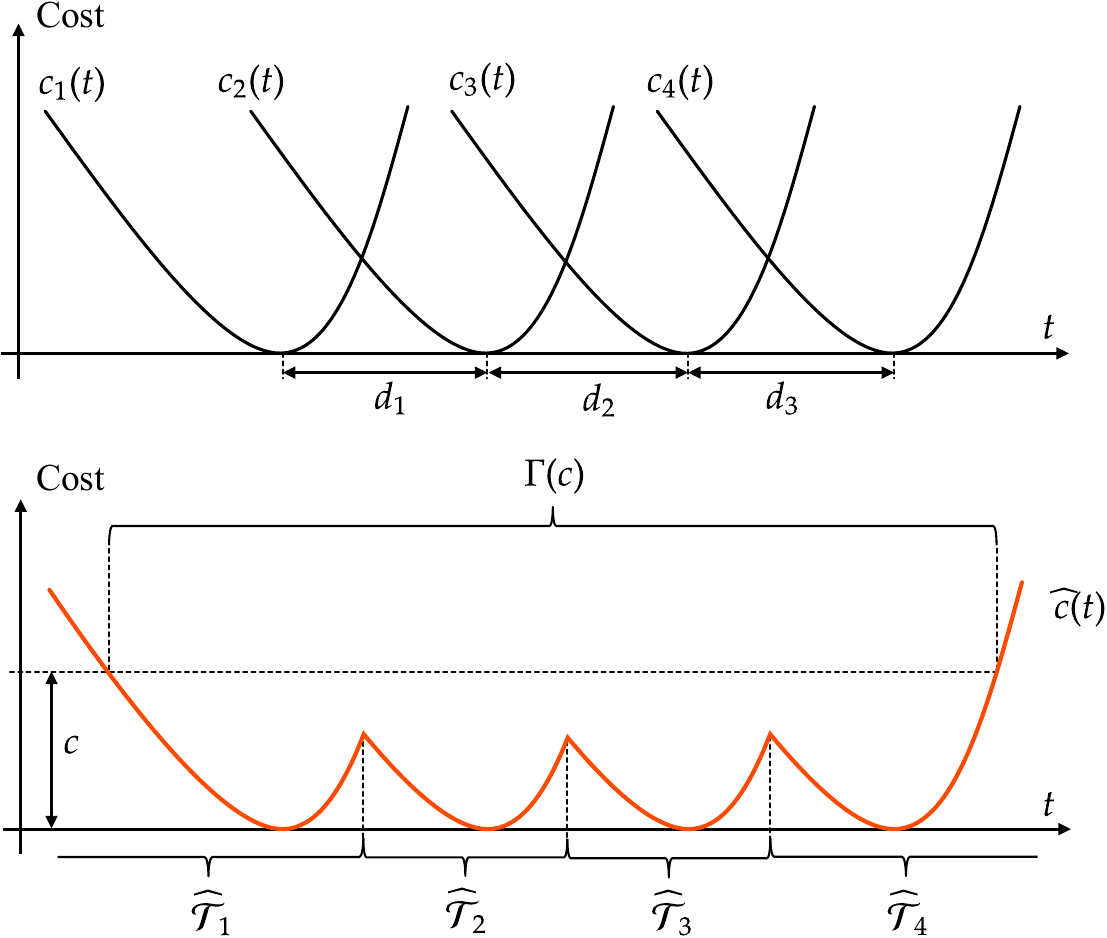}
  \caption{$\widehat{c}(t)$, $\widehat{\ClT}_{k}$, and $\Gamma(c)$ in Relation to the Schedule Delay Cost Function $c_{k}(t)$.}
  \label{fig:SDCF}
\end{figure}

The workers live in the residential location and engage with the firms located in the CBD.
The workers live in the residential location and engage with the firms located in the CBD.
We assume that the land area consumed by all workers is a fixed unit size (i.e., lot size is fixed), and there are no crowding or other externalities from residential population density.
The workers choose the residential location $i \in \ClI$, the office work ratio $h \in \ClH$, preferred arrival time $t_{k}$, $k \in \ClK$ and actual arrival time $t \in \ClT$ at the CBD to maximize their utility, where $\ClT$ is the set of all possible arrival times (assignment period).
\footnote{
  By imposing restrictions on each worker's choice of the preferred arrival time $t_{k}$ (e.g., a worker can only choose $t_{1}$ or $t_{2}$), we can model the firm's working hours requirements (i.e., core working hours).
  In order to simplify the analysis, this study assumes that workers are free to choose, from the set of preferred arrival times, their own preferred arrival time.}
  Hereafter, we call the $(h,i,k,t)$-worker as the worker who lives in the residential location $i$, chooses the preferred arrival time $t_{k}$ and actual arrival time $t$ for convenience.
  The utility is calculated by the wage, commuting cost and land rent as the sum of $N$ days, where $N$ is the number of working days per unit term.
  The land rent at the residential location $i \in \ClI$ is denoted by $r_{i}$, which is determined by the land market (market clearing condition is given later) and paid from the workers to an absentee landlord.
  We assume that there is one representative absentee landlord.  
The commuting cost includes the free-flow travel, queueing delay and the schedule delay costs.
  The queuing delay at bottleneck $i$ for the commuter arriving at time $t$ at the destination is denoted by $w_{i}(t)$.
The commuting cost for an individual worker whose preferred arrival time is $t_{k}$, actual arrival time is $t$, residential location is $i$ is as follows:
\begin{align}
  &C_{i,k}(t) = c_{k}(t) + \alpha \sum_{j;j \leq i} \left( w_{j}(t) + f_{j} \right)
  &&\forall k \in \ClK,
  \quad \forall i \in \ClI, \quad \forall t \in \ClT,
  \label{eq:commuting_cost}
\end{align}
where $\alpha$ is a parameter representing the value of time.
We assume $\alpha=1$ for all commuters.
The function $c_{k}(t)$ is the schedule delay cost function.
We assume that $c_{k}(t)$ can be expressed as $c(t-t_{k})$, where ${c}(T)$ denotes a base schedule delay cost function with $c(0)=0$.
This function $c(T)$ is assumed to be convex and piecewise differentiable with respect to $T$.
Moreover, we introduce several functions and sets related to $c_{k}(t)$ as follows.
\begin{align}
    &\widehat{c}(t) \equiv \inf_{k \in \ClK} \left\{ c_{k}(t) \right\},
    \label{eq:widehat_c}
    \\
    &\widehat{\ClT}_{k} \equiv
    \left\{ t \in \ClT \mid k \in \argmin_{k \in \ClK} c_{k}(t) \right\},
    \label{eq:widehat_T}
    \\
    &\Gamma(c) \equiv \left\{ t \in \ClT \mid \widehat{c}(t) \leq c  \right\}.
    \label{eq:gamma_c}
\end{align}
The function $\widehat{c}(t)$ denotes the envelope function of $\{ c_{k}(t) \}_{k \in \ClK}$.
The set $\widehat{\ClT}_{k}$ denotes the set of times at which the function $c_{k}(t)$ attains its minimum value among $\{ c_{k}(t) \}_{k \in \ClK}$.
The set $\Gamma(c)$ denotes the set of times for which the envelope function $\widehat{c}(t)$ is less than or equal to $c$.
The relationships between these functions and sets are depicted in \Cref{fig:SDCF}.
In addition, we introduce the following function $\overline{c}(X, \mu)$:
\begin{align}
  &\overline{c}(X, \mu) \equiv c
  &&\mbox{where } c \mbox{ is the solution to} \quad \int_{t \in \Gamma(c)} \mu \mathrm{d} t = X.
  \label{eq:overline_c}
\end{align}
The function $\overline{c}(X, \mu)$ is used when we derive the optimal and equilibrium solutions to the short-term model.
The value $c$ of this function represents the schedule delay cost corresponding to a time window $\Gamma(c)$, as shown in \Cref{fig:SDCF}.
The time window $\Gamma(c)$ for any demand $X$ and capacity $\mu$ is determined by $\int_{t \in \Gamma(c)} \mu \mathrm{d} t = X$, which means that the maximum demand of commuters who can pass through the capacity $\mu$ bottleneck with the schedule cost less than or equal to $c$ is $X$.
\par
The utility function for worker is given by:
\begin{align}
  N \times U_{i,k}(h, t) 
  &\equiv
  N \left\{ h \left( \theta^{O} - C_{i,k}(t)\right)
  + \left( 1 - h \right) \theta^{R}
  - r_{i} \right\}
  \notag
  \\
  &= 
  N \Biggl\lbrace h \left( \theta^{O} - c_{k}(t) - \sum_{j;j\leq i} 
  \left( w_{j}(t) + f_{j}  \right) \right) 
  + \left( 1 - h \right) \theta^{R} - r_{i}  \Biggr\rbrace.
  \label{eq:utility_function}
\end{align}
The first term on the right-hand side of~\eqref{eq:utility_function} denotes the utility associated with office work days, the second term represents the utility  associated with remote work days and the third term correspond to the land rent.
\par
Let $Q_{i}(h)$ denote the number of workers with an office work ratio of $h$ and residing in location $i$, that is, $Q = \sum_{i \in \ClI} \int_{h \in \ClH} Q_{i}(h) \mathrm{d} h$.
Furthermore, let $q_{i,k}(h, t)$ denote the inflow rate of workers with an office work ratio of $h$, residing in location $i$, and having a preferred arrival time of $k$ and actual arrival time of $t$.
Both $\{Q_{i}(h)\}$ and $\{q_{i,k}(h, t)\}$ are endogenous variables determined by the equilibrium condition.

\subsection{Short-term Model}
The short-term equilibrium model represents the within-day dynamic traffic flow, which consists of the commuter conservation condition, queueing condition, and preferred and actual arrival times choice condition.
In order to clearly formulate the short-term equilibrium, we employ the \textit{Lagrangian-like coordinate system approach}~\citep{Kuwahara1993-vt,Akamatsu2015-ip,Akamatsu2021-zg}.
In this approach, variables are expressed regarding the arrival time at the destination rather than the origin or bottleneck.
This representation is particularly suitable for considering the ex-post travel time of each commuter during their trip, enabling the straightforward tracing of time-space paths.
\par
Several variables are introduced to describe the within-day dynamic traffic flow.
We define $\tau_{i}(t)$ and $\sigma_{i}(t)$ as the arrival and departure times, respectively, at bottleneck $i\in\ClI$ for commuters whose destination arrival time is $t$.
These variables satisfy the following relationships:
\begin{align}
  &\tau_{i}(t) = t - \sum_{j; j \leq i} \left( w_{j}(t) + f_{j} \right)
  &&\forall i \in \ClI, \quad \forall t \in \ClT,
  \label{eq:tau=s-sumw}
  \\
  &\sigma_{i}(t) = t - \sum_{j; j < i}  \left( w_{j}(t) + f_{j} \right)
  &&\forall i \in \ClI, \quad \forall t \in \ClT.
  \label{eq:sigma=s-sumw}
\end{align}  
Note that $\tau_{i}(t)$ must satisfy the following relationship:
\begin{align}
  &\cfrac{\mathrm{d}\tau_{i}(t)}{\mathrm{d}t}\equiv \dot{\tau}_{i}(t) > 0,
  &&\forall i \in \ClN, \quad \forall t \in \ClT,
  \label{eq:TauMonotone}
\end{align}
where an overdot represents the derivative of the variable with respect to the destination arrival time $t$.
This condition ensures the Lipschitz continuity of cumulative arrival flow, ensuring the physical consistency of traffic flow~\citep{Akamatsu2015-ip,Fu2022-nl,Sakai2022-vm}.
Henceforth, we will refer to this condition as the consistency condition.
\par
  Building upon the Lagrangian-like coordinate system approach, we express the commuter conservation condition, queuing condition, and choice conditions for preferred and actual arrival times as follows:
\begin{align}
&Q_{i}(h)=  \sum_{k \in \ClK} \int_{t \in \ClT} q_{i,k}(h,t) \mathrm{d}t
&&
\forall i \in \ClI, \quad \forall h \in \ClH,
\label{eq:DUE_CommuterCnsv}
\\
&\begin{dcases}
  \sum_{j;j \geq i} \sum_{k \in \ClK} \int_{h \in \ClH}  h q_{j,k}(h,t) \mathrm{d}h
  = \mu_{i} \dot{\sigma}_{i}(t)
  &\mathrm{if}\quad 
  w_{i}(t) > 0
  \\
  \sum_{j;j \geq i} \sum_{k \in \ClK} \int_{h \in \ClH}  h q_{j,k}(h,t) \mathrm{d}h
  \leq \mu_{i} \dot{\sigma}_{i}(t)
  &\mathrm{if}\quad 
  w_{i}(t) = 0
 \end{dcases}
 &&
 \forall i \in \ClN,
 \quad \forall t \in \ClT,
 \label{eq:DUE_QueueingCondition}
  \\
  &\begin{dcases}
  c_{k}(t) + \sum_{j;j\leq i} w_{j}(t)  =  \lambda_{i}(h)
  &\mathrm{if}\quad  q_{i,k}(h,t) > 0
  \\
  c_{k}(t) + \sum_{j;j\leq i} w_{j}(t)  \geq \lambda_{i}(h)
  &\mathrm{if}\quad  q_{i,k}(h,t) =  0
  \end{dcases}
  &&\forall k \in \ClK, 
  \quad \forall i \in \ClI,
  \quad \forall t \in \ClT,
  \quad \forall h \in \ClH,
  \label{eq:DUE_DeparturePrefferd_Time_Choice}
\end{align}
where $\lambda_{i}(h)$ represents the equilibrium commuting cost, excluding the free-flow travel cost for workers with an office work ratio of $h$ and residing in location $i$.
The details of the queueing condition\footnote{
  We assume that the office work schedules for all workers are leveled.
  That is, for all days, the number of commuters to the CBD whose office work ratio is $h$ and the residential location is $i$ is $h Q_{i,h}$.
  Therefore, the arrival flow of commuters to the CBD whose office work ratio is $h$ and the residential location is $i$ is $\sum_{k \in \ClK} \int_{h \in \ClH} h q_{i,k}(h,t) \mathrm{d}t$.
} based on the Lagrangian-like coordinate system approach are given by \citet{Akamatsu2015-ip,Fu2022-nl,Sakai2022-vm}.
\par
We formally define the short-term equilibrium as follows:
\begin{definition}(Short-term equilibrium)
  \label{dfn:short-term_equilibrium}
  The short-term equilibrium is a collection of variables $(
   \UE{\Vtq} \equiv \{ \UE{q}_{i,k}(h,t) \}$, 
   $\UE{\Vtw} \equiv \{ \UE{w}_{i}(t) \}$, 
  $\UE{\Vtlambda}  \equiv \{ \UE{\lambda}_{i}(h) \})$ that satisfies \eqref{eq:DUE_CommuterCnsv}, \eqref{eq:DUE_QueueingCondition}, and \eqref{eq:DUE_DeparturePrefferd_Time_Choice}.
\end{definition}
\par
Additionally, we introduce the \textit{short-term optimal problem}, which aims to find a socially optimal traffic flow pattern.
The decision variable in this problem is $\Vtq = \{ q_{i,k}(h,t) \}$, and the objective function is the total travel cost. The problem is formulated as follows:
\begin{align}
  \text{[ST-SO]} &\quad 
  \notag
  \\
  \min_{\Vtq \geq \Vt0}. \quad 
  &Z( \Vtq \mid \VtQ)  =  \sum_{i \in \ClI}
  \sum_{k \in \ClK} \int_{t \in \ClT} \int_{h \in \ClH}
  h  c_{k}(t)  q_{i,k}(h,t) \mathrm{d}h \mathrm{d} t
  \\
  \text{s.t.} \quad 
  &\sum_{k \in \ClK} \int_{t \in \ClT} q_{i,k}(h,t) 
  \mathrm{d} t = Q_{i}(h)
  && 
  \forall i \in \ClI, \quad \forall h \in \ClH,
  \label{eq:DSO_DemandCnsv}
  \\
  &\sum_{j;j\geq i} \sum_{k \in \ClK}  \int_{h \in \ClH} 
  h q_{i,k}(h,t) \mathrm{d}h \leq \mu_{i}
  && 
  \forall i \in \ClI, \quad \forall t \in \ClT.
  \label{eq:DSO_CapacityCnst}
\end{align}
Constraint \eqref{eq:DSO_DemandCnsv} represents the workers conservation condition, and constraint \eqref{eq:DSO_CapacityCnst} represents the capacity constraint to prevent bottleneck queues.
The optimality condition for this problem is given by:
\begin{align}
&Q_{i}(h) = \sum_{k \in \ClK} \int_{t \in \ClT}  \SO{q}_{i,k}(h,t) \mathrm{d} t
&&\forall i \in \ClI, \quad \forall h \in \ClH,
\label{eq:DSO_CommuterCnsv}
\\
&\begin{dcases}
  \sum_{j;j\geq i} \sum_{k \in \ClK} \int_{h \in \ClH} h \SO{q}_{j,k}(h,t) \mathrm{d}h=  \mu_{i}
  &\mathrm{if}\quad  \SO{p}_{i}(t) > 0
  \\
  \sum_{j;j\geq i} \sum_{k \in \ClK} \int_{h \in \ClH} h \SO{q}_{j,k}(h,t) \mathrm{d}h \leq \mu_{i} 
  &\mathrm{if}\quad  \SO{p}_{i}(t) = 0
\end{dcases}
&&
\forall i \in \ClI, \quad \forall t \in \ClT,
\label{eq:DSO_OptimalTolling}
\end{align}
\begin{align}
&\begin{dcases}
   c_{k}(t) + \sum_{j;j\leq i} \SO{p}_{j}(t)  =  \SO{\lambda}_{i}(h) 
&\mathrm{if}\quad  \SO{q}_{i,k}(h,t) > 0
\\
c_{k}(t) + \sum_{j;j\leq i} \SO{p}_{j}(t)  \geq \SO{\lambda}_{i}(h) 
&\mathrm{if} \quad \SO{q}_{i,k}(h,t) =  0
\end{dcases}
&&
\forall k \in \ClK, 
\quad \forall i \in \ClI,
\quad \forall t \in \ClT,
\quad \forall h \in \ClH,
\label{eq:DSO_DeparturePrefferd_Time_Choice}
\end{align}
  where $\lambda_{i}(h)$ and $p_{i}(t)$ are the Lagrange multipliers corresponding to constraints \eqref{eq:DSO_DemandCnsv} and \eqref{eq:DSO_CapacityCnst}, respectively.
  As mentioned in previous studies~\citep{Arnott1990-ta,Laih1994-hi,Lindsey2012-gg,Chen2015-ku,Osawa2018-hg,Fu2022-nl,Sakai2022-ay}, 
  these optimality conditions can be interpreted as equilibrium conditions under an optimal dynamic congestion pricing scheme.
  Under this interpretation, the Lagrange multiplier $\lambda_{i}(h)$ represents the equilibrium commuting cost for workers residing in location $i$ and opting an office work ratio of $h$. 
  The Lagrange multiplier $p_{i}(t)$ represents the optimal pricing pattern at bottleneck $i$ at time $t$.
\footnote{
  Another interpretation of $\{ \SO{\Vtp}(t) \}_{t \in \ClT}$ is the market clearing price pattern under a time-dependent tradable bottleneck permit scheme~\citep{Wada2013-li,Akamatsu2017-bi}. 
}
\par
We formally define the short-term optimum as follows:
\begin{definition}(Short-term optimum)
  \label{dfn:short-term_optimum}
  The short-term optimum is a collection of variables $(
    \SO{\Vtq} \equiv \{ \SO{q}_{i,k}(h,t) \}$, 
    $\SO{\Vtp} \equiv \{ \SO{p}_{i}(t) \}$, 
    $\SO{\Vtlambda}  \equiv \{ \SO{\lambda}_{i}(h) \})$
    that is optimal solution and optimal Lagrangian multiplier to [ST-SO]. 
\end{definition}
  In this paper, depending on the context, the short-term optimal state defined by \Cref{dfn:short-term_optimum} may be referred to as the \textit{toll equilibrium}, and the short-term equilibrium state defined by \Cref{dfn:short-term_equilibrium} may be referred to as the \textit{no-toll equilibrium}.

\subsection{Long-term Model}
  The long-term equilibrium problem aims to determine the equilibrium location and office work ratio pattern, incorporating the worker conservation condition, land market clearing condition, and location and office work ratio choice condition. We assume a standard, perfectly competitive land market at each location.
  We express the long-term equilibrium conditions as follows:
\begin{align}
  &\sum_{i \in \ClI} \int_{h \in \ClH} Q_{i}(h)
  \mathrm{d}h = Q,
  \label{eq:DUE_WorkerCnsv}
  \\
  &\begin{dcases}
    \int_{h \in \ClH} Q_{i}(h)
    \mathrm{d}h = A_{i}
    &\mathrm{if}\quad  r_{i} > 0
    \\
    \int_{h \in \ClH} Q_{i}(h)
    \mathrm{d}h \leq A_{i}
    &\mathrm{if}\quad  r_{i} = 0
  \end{dcases}
  &&\forall i \in \ClI,
  \label{eq:DUE_LandMarket}
  \\
  &
  \begin{dcases}
    h  \left( \theta^{O} - \lambda_{i}(h)- \sum_{j;j\leq i} f_{j} \right)
    + \left( 1 - h  \right) \theta^{R} - r_{i}= \rho 
  &\mathrm{if}\quad  Q_{i}(h)> 0
  \\
  h  \left( \theta^{O} - \lambda_{i}(h)- \sum_{j;j\leq i} f_{j} \right)
  + \left( 1 - h  \right) \theta^{R} - r_{i}\leq \rho 
  &\mathrm{if}\quad  Q_{i}(h)= 0
  \end{dcases}
  &&
  \forall i \in \ClI,
  \quad \forall h \in \ClH,
  \label{eq:DUE_WorkingstyleLocation_Choice}
\end{align}
where $\rho$ is the equilibrium utility of workers.
\par
We formally define the long-term equilibrium as follows:
\begin{definition}(Long-term equilibrium)
  \label{dfn:long-term_equilibrium}
  A long-term equilibrium is a collection of variables $(
    \UE{\VtQ} \equiv \{ \UE{Q}_{i}(h) \}$,
    $\UE{\Vtr} \equiv \{ \UE{r}_{i} \}$,
    $\UE{\rho}$ that satisfies \eqref{eq:DUE_WorkerCnsv}, \eqref{eq:DUE_LandMarket}, and \eqref{eq:DUE_WorkingstyleLocation_Choice}.
\end{definition}

\subsection{Integrated equilibrium}
We introduce the integrated equilibrium, which combines the short- and long-term equilibrium problems, as follows:
\begin{definition}[Integrated equilibrium]
  \label{dfn:integrated_equilibrium}
  An integrated equilibrium is a collection of variables
  $(
    \UE{\Vtq} \equiv \{ \UE{q}_{i,k}(h,t) \}$, 
    $\UE{\Vtp} \equiv \{ \UE{p}_{i}(t) \}$, 
    $\UE{\Vtlambda}  \equiv \{ \UE{\lambda}_{i}(h) \},
    \UE{\VtQ} \equiv \{ \UE{Q}_{i}(h) \},
    \UE{\Vtr} \equiv \{ \UE{r}_{i} \},
    \UE{\rho}
  )$ that satisfies 
  \eqref{eq:DUE_CommuterCnsv}, 
  \eqref{eq:DUE_QueueingCondition},
  \eqref{eq:DUE_DeparturePrefferd_Time_Choice}, 
  \eqref{eq:DUE_WorkerCnsv},
  \eqref{eq:DUE_LandMarket}, 
  and \eqref{eq:DUE_WorkingstyleLocation_Choice}.
\end{definition}
\Cref{fig:Strucure_IntModel} illustrates the structure of the integrated model, which encompasses both short- and long-term models. 
These models are interconnected through the variables $\{ Q_{i}(h) \}$ and $\{\lambda_{i}(h)\}$.
The short-term model is formulated given equilibrium commuting demand patterns $\{ Q_{i}(h) \}$ determined by the long-term equilibrium.
Conversely, the long-term model is formulated given the equilibrium commuting cost pattern $\{\lambda_{i}(h)\}$ determined by the short-term equilibrium.
\par
This integrated equilibrium model can also describe equilibrium states without TLC or SWH. 
When the set of office work ratios is reduced to $\ClH = \{1\}$, our model describes the equilibrium state without TLC.
Similarly, when the size of the set of preferred arrival times is reduced to $K=1$, our model describes the equilibrium state without SWH. 
In light of this, in the remainder of the paper, we first derive a general solution to the integrated equilibrium model; then, we compare these schemes' long- and short-term effects using the general solution.

\section{Short-term Equilibrium}
\label{sec:ShortTermEquilibrium}
In this section, we present the closed-form solution to the short-term problem. 
In \Cref{subsec:QRP}, we outline the approach to solving the problem and introduce the queue replacement principle (QRP), which plays a pivotal role in addressing short-term problem. 
The QRP represents the theoretical relationship between equilibrium and optimum, and it enables us to derive the short-term \textit{equilibrium} state by solving the short-term \textit{optimum} problem.
Therefore, in \Cref{subsec:ST-optimum}, we first derive the solution to the short-term optimal problem.
Subsequently, by using the QRP and the solution to the short-term optimal problem, we construct the solution to the short-term equilibrium in \Cref{subsec:ST-equilibrium}.

\subsection{Queue Replacement Principle Approach}
\label{subsec:QRP}
In order to solve the short-term equilibrium problem, we use the QRP approach for the dynamic user equilibrium problem in the corridor network~\citep{Fu2022-nl,Sakai2022-ay,Sakai2023-gs}.
The QRP enables us to derive the equilibrium (no toll equilibrium) state by using the optimal (toll equilibrium) state.
We formally introduce the QRP as follows:
\begin{definition}
  (QRP; \citealp{Fu2022-nl,Sakai2022-ay,Sakai2023-gs})
  \label{dfn:QRP_Multiple}
  If there exists a no toll equilibrium state in which the queueing delay pattern is equal to the optimal pricing pattern:
  \begin{align}
    \UE{w}_{i}(t) = \SO{p}_{i}(t)
    &&\forall i \in \ClN, \quad \forall t \in \ClT,
  \end{align}
  then, the QRP holds.
\end{definition}
\noindent
This relationship implies two important facts.
First, the optimal state can be achieved by imposing a congestion toll equal to the queueing pattern in the no-toll equilibrium state.
Second, the congestion toll required to attain the optimal state is equal to the queuing pattern in the no-toll equilibrium state. Despite these facts appearing similar at first glance, the latter view is more mathematically sound. 
This is because solving the optimization problem is generally easier than directly solving complex equilibrium problems. When the QRP holds, we can construct the no-toll equilibrium state using the toll equilibrium (optimal) state without directly solving complex equilibrium problems.
\par
Therefore, we first analytically derive the short-term optimal state, as defined by \Cref{dfn:short-term_optimum}.
We then prove that the QRP holds under an assumption.
Finally, we construct the short-term equilibrium state, as defined by \Cref{dfn:short-term_equilibrium}, based on the solution to the short-term optimal problem using the QRP.

\subsection{Solution to the Short-term Optimum Problem}
\label{subsec:ST-optimum}
  In this section, we present the closed-form solution to the short-term optimal problem [ST-SO].
  As discussed in the next section, the workers who live far from the CBD do not choose the office work with commuting to the CBD in the long-term equilibrium.
  Therefore, the analysis of the short-term equilibrium, which represents the commuting problem, focuses on workers residing in the CBD and its vicinity.
  The identification of the boundary point will be addressed in the analysis of the long-term equilibrium in the subsequent section, and the point is denoted here as $i^{\ast}$.
  We define $\ClI^{\ast} = \{ i^{\ast}, i^{\ast}-1, \ldots, 1 \}$ as the set of residential locations in which the workers choose the office work in the long-term equilibrium.
  The analysis of the short-term equilibrium in this section is limited to the residential locations in $\ClI^{\ast}$.
\par
In order to derive the analytical solution, we first decompose the short-term optimal problem into bottleneck-based problems under a specific assumption.
These bottleneck-based problems share the same structure as the analytically solvable single bottleneck problem. 
Based on this, we derive closed-form solutions for each decomposed problem independently. 
Finally, we construct the closed-form solution for the short-term optimal problem by combining the solutions to bottleneck-based problems.
\par
In order to ensure clarity and exclude arbitrary optimal flow patterns in the corridor problem, we introduce the following assumptions regarding flow and pricing patterns in the optimal solution:
\begin{assumption}
	For all $i \in \ClI^{\ast}$ and all $t \in \ClT$, the following relationship is satisfied:
	\begin{align}
		\SO{p}_{i}(t) > 0
		\quad \Leftrightarrow \quad 
		 \sum_{k \in \ClK} \int_{h \in \ClH} \SO{q}_{i,k}(h,t) \mathrm{d} h  > 0.
		\label{eq:Asm_p>0q>0}
	\end{align}
  The time window $\SO{\ClT}_{i} \equiv \mathrm{supp} \left( \SO{p}_{i}(t) \right) = \mathrm{supp} \left( \SO{q}_{i}(t) \right)$ is convex for all $i \in \ClI^{\ast}$.
	\label{asm:p>0q>0}
\end{assumption}
\noindent
This assumption implies that all commuters experience non-zero optimal congestion prices at each bottleneck along their route, ruling out false bottlenecks in which optimal prices are zero at any time. 
In a corridor network with false bottlenecks, the optimal flow pattern becomes arbitrary~\citep{Fu2022-nl}. 
By considering this assumption, we can exclude cases where the closed-form solutions would entail arbitrary optimal flow patterns, thus simplifying the analysis.
\par
Under this assumption, we observe the following property regarding the arrival time windows of commuters departing from bottleneck $i$ and $i+1$:
	\begin{align}
		&\SO{\ClT}_{i} \subset \SO{\ClT}_{(i+1)}
		&&\forall i \in \ClI^{\ast} \setminus \{ I \}.
		\label{eq:Lem_S_i_subset_S_i+1}
	\end{align}
Using this property, we can determine the aggregated short-term optimal flow pattern, as follows: 
\begin{lemma}
	Suppose that \Cref{asm:p>0q>0} holds, then the aggregated short-term optimal flow is given by:
	\begin{align}
	&\sum_{k \in \ClK} \int_{h \in \ClH} h q_{i,k}(h,t)
    \mathrm{d}h
		\leq \overline{\mu}_{i}
	 &&\forall i \in \ClI^{\ast}, \quad 
   \forall t \in \ClT,
	 \label{eq:Lem_q=mu}
	\end{align}
    where $\overline{\mu}_{i} = \mu_{i} - \mu_{(i+1)}$ and $\mu_{i^{\ast} + 1} = 0$.
	\label{lem:q=mu}
\end{lemma}
\noindent{\bf Proof}.
  See Appendix \ref{subsec:proof_lem_q=mu}.
\par
This lemma implies that the capacity constraints of each bottleneck can be handled independently, enabling us to decompose [ST-SO] into bottleneck-based problems. 
By solving the following bottleneck-based problem, we can obtain the short-term optimal solution.
\begin{align}
  \min_{ \{q_{i,k}(h,t) \}_{k,h,t} } \quad &
  \sum_{k \in \ClK} \int_{h \in \ClH} \int_{t \in \ClT} h c_{k}(t) q_{i,k}(h,t) \mathrm{d}t \mathrm{d}h
  \\
  \text{s.t.} \quad 
  &\sum_{k \in \ClK} \int_{t \in \ClT} q_{i,k}(h,t) 
  \mathrm{d} t = Q_{i}(h)
  &&\forall h \in \ClH
  \quad [\lambda_{i}(h)],
  \\
  &\sum_{k \in \ClK} \int_{h \in \ClH} 
  h q_{i,k}(h,t) \mathrm{d}h
  \leq \mu_{i} - \mu_{i+1}
  &&\forall t \in \ClT
  \quad [P_{i}(t)].
\end{align}
This bottleneck-based problem can be solved as a single bottleneck problem, which can be analytically solved.
  This is summarized in the following lemma.
\begin{lemma}
  The variable set $\{ q^{\ast}_{i,k}(h, t) \}$, $\{ \lambda^{\ast}(h)(t) \}$ and $\{ P^{\ast}_{i}(t) \}$ is the solution to the bottleneck-based problem.
  \begin{align}
    &\lambda^{\ast}_{i} = \overline{c}(X^{\ast}_{i}, \overline{\mu}_{i}), 
    \\
    &\lambda^{\ast}_{i}(h) = h \lambda^{\ast}_{i}
    &&\forall h \in \ClH,
    \\
    &P^{\ast}_{i}(t)
     = \max \left\{ 0, \lambda^{\ast}_{i} - \widehat{c}(t) \right\}
    &&\forall t \in \ClT,
    \\
    &\int_{h \in \ClH} hq_{i,k}(h,t) 
    \mathrm{d}h = 
    \begin{dcases}
      \overline{\mu}_{i}
      &\mathrm{if}\quad  
      t \in \SO{\ClT}_{i} \cap \widehat{\ClT}_{k}
      \\
      0
      \quad &\mathrm{otherwise}
    \end{dcases}
    && 
    \forall k \in \ClK, \forall t \in \ClT,
    \\
    &\mbox{where } X^{\ast}_{i} = \int_{h \in \ClH} hQ_{i}(h) \mathrm{d}h.
  \end{align}
  Here, $\overline{c}(X, \mu)$, $\widehat{c}(t)$, and $\widehat{\ClT}_{k}$ are defined in \eqref{eq:overline_c}, \eqref{eq:widehat_c}, and \eqref{eq:widehat_T}, respectively.
\label{lem:Solution_BottleneckBasedProblem}
\end{lemma}
\noindent{\bf Proof}.
  See Appendix \ref{subsec:proof_lem_Solution_BottleneckBasedProblem}.
\par
\begin{figure}[tbp]
  \center
     \includegraphics[clip, width=0.7\columnwidth]{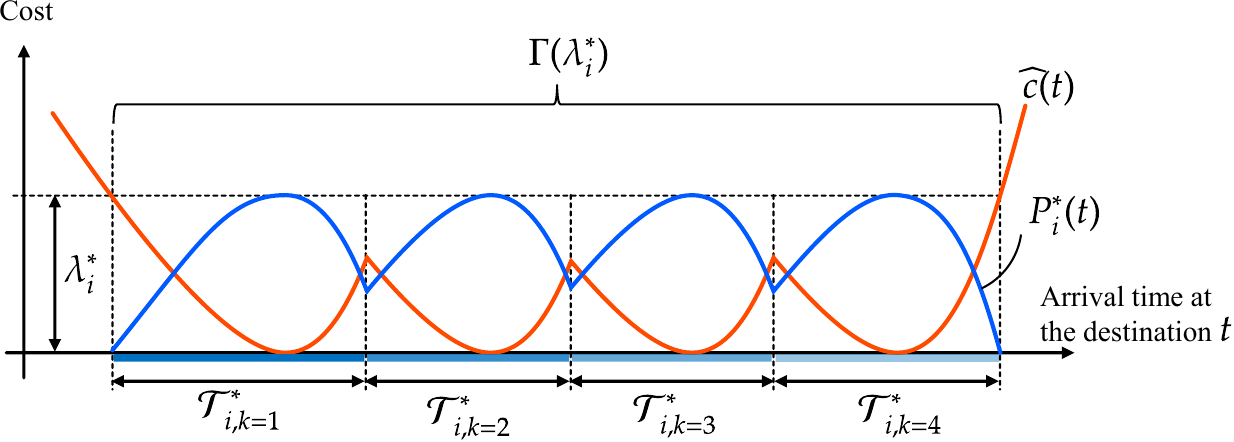}
  \caption{Solution to the Bottleneck-based Problem.}
  \label{fig:SolutoinToSingle}
\end{figure}

  The optimal solution $\Vtq$ is not unique; however, the optimal solution $\Vtlambda$ and $\VtP$ are unique.
  \Cref{fig:SolutoinToSingle} illustrates the solution to the bottleneck based problem.
  Let $\ClT^{\ast}_{i,k}$ be the set of times where the optimal flow $q^{\ast}_{i,k}(h,t)$ is positive, i.e., $\ClT^{\ast}_{i,k} = \SO{\ClT}_{i} \cap \widehat{\ClT}_{k}$.
  This solution implies that first-in-first-work (FIFW) principle \citep{Daganzo1985-ls} holds for the optimal flow.
  Moreover, $P^{\ast}_{i}(t)$ is determined by the difference between $\lambda^{\ast}_{i}$ and $\widehat{c}(t)$ for all $t \in \SO{\ClT}_{i}$.

The entire closed-form solution can be derived by combining the analytical solutions for the decomposed problems. 
The flow variable $\Vtq$ and the equilibrium commuting cost $\Vtlambda$ remain the same as in the bottleneck-based problem, while the optimal congestion price needs to be calculated separately, as follows:
\begin{align}
  &\SO{p}_{i}(t)
  = P^{\ast}_{i}(t) - P^{\ast}_{i-1}(t)
  &&\forall t \in \ClT.
\end{align}
This is summarized in the following lemma.
\begin{lemma}(Solution to the short-term optimal problem)
  \label{lem:Solution_ST-SO}
  The following variables are the solution to the problem [ST-SO]:
  \begin{align}
    &\SO{\lambda}_{i} = \overline{c}(X^{\ast}_{i}, \overline{\mu}_{i}), 
    &&\forall i \in \ClI,
    \\
    &\SO{\lambda}_{i}(h) = h \SO{\lambda}_{i}
    &&\forall i \in \ClI, \quad \forall h \in \ClH,
    \\
    &\SO{p}_{i}(t)
     = \max \left\{ 0, \lambda^{\ast}_{i} - \widehat{c}(t) 
     - \sum_{j;j>i} \SO{p}_{j}(t) \right\}
     &&\forall i \in \ClI, \quad \forall t \in \ClT,
    \\
    &\int_{h \in \ClH} h\SO{q}_{i,k}(h,t) \mathrm{d}h= 
    \begin{dcases}
      \overline{\mu}_{i}
      &\mathrm{if}\quad  
      t \in \SO{\ClT}_{i} \cap \widehat{\ClT}_{k}
      \\
      0
      \quad &\mathrm{otherwise}
    \end{dcases}
    &&\forall k \in \ClK, \quad \forall i \in \ClI, \quad \forall t \in \ClT,
    \\
    &\mbox{where } X^{\ast}_{i} = \int_{h \in \ClH} hQ_{i}(h) \mathrm{d}h.
  \end{align}
\end{lemma}

\subsection{Solution to the Short-term Equilibrium Problem}
\label{subsec:ST-equilibrium}
  By using the obtained short-term optimum, we now derive the closed-form equilibrium solution.
  Specifically, we substitute the derived optimal pricing pattern into the short-term equilibrium conditions as a queuing delay pattern. 
  We demonstrate that this queuing delay pattern satisfies all short-term equilibrium conditions without contradiction under certain conditions. 
  This approach simultaneously validates the QRP and derives the equilibrium solution using the optimal solution.
\par
In conclusion, we find that the QRP holds under an assumption of the schedule delay cost function.
Thus, using the short-term optimal solution shown in \Cref{lem:Solution_ST-SO}, we can obtain the solution to the short-term equilibrium problem under the assumption.
This is summarized in the following propositions. 
\begin{proposition}
  Suppose that the schedule delay cost function satisfies the following conditions:
  \begin{align}
    &-1 < \dot{c}_{k}(t) < \dfrac{\mu_{i}-\mu_{i+1}}{\mu_{i+1}}
    &&\forall k \in \ClK, \quad \forall i \in \ClN, \quad \forall t \in \ClT.
    \label{eq:QRP_Condition}
  \end{align}
  Then, the optimal pricing pattern in the short-term optimal state is equal to the queueing delay pattern in the short-term equilibrium state:
  \begin{align}
    &\SO{p}_{i}(t) = \UE{w}_{i}(t)
    &&\forall i \in \ClI, \quad \forall t \in \ClT.
  \end{align}
  \label{pro:QRP}
\end{proposition}
\noindent{\bf Proof}. See Appendix \ref{subsec:prf_QRP_DUE_Solution}
\begin{proposition}(Solution to the short-term equilibrium problem)
  \label{pro:DUE_Solution}
  Suppose that the schedule delay cost function satisfies \eqref{eq:QRP_Condition}.
  Then, the following variables $\{\UE{\lambda}_{i}(h)\}$, $\{\UE{w}_{i}(t)\}$ and $\{\UE{q}_{h,i,k}(t)\}$ are the solution to the short-term equilibrium problem.
  \begin{align}
    &\UE{w}_{i}(t) = \SO{p}_{i}(t)
    &&\forall i \in \ClI, \quad \forall t \in \ClT,
    \\
    &\UE{\lambda}_{i}(h) = \SO{\lambda}_{i}(h)
    &&\forall h \in \ClH, \quad \forall i \in \ClI,
    \\
    &\UE{\ClT}_{i} = \SO{\ClT}_{i}
    &&\forall i \in \ClI,
    \\
    &\int_{h \in \ClH} h \UE{q}_{i,k}(h,t) 
    \mathrm{d}h 
    = 
    \begin{dcases}
        \mu_{i} \UE{\dot{\sigma}}_{i}(t)
         - 
        \mu_{i+1} \UE{\dot{\sigma}}_{i+1}(t)
      &\mathrm{if}\quad  t \in \UE{\ClT}_{i} \cap \widehat{\ClT}_{k}
      \\
      0
      \quad &\mathrm{otherwise}
    \end{dcases}
    &&\forall k \in \ClK, \quad \forall i \in \ClI,
    \label{eq:DUE_solution_q}
  \end{align}
  where $\{\UE{q}_{i,k}(h,t)\}$ satisfies the following conditions:
  \begin{align}
    &\sum_{k \in \ClK} \int_{t \in \ClT} \UE{q}_{i,k}(h,t) \mathrm{d}t = Q_{i}(h)
    &&\forall i \in \ClI, \quad  \forall h \in \ClH.
    \label{eq:DUE_solution_q_sum}
  \end{align}
\end{proposition}
\noindent{\bf Proof}. See Appendix \ref{subsec:prf_QRP_DUE_Solution}
\\
  The detail technical reason is given in Appendix \ref{subsec:prf_QRP_DUE_Solution}, the condition \eqref{eq:QRP_Condition} prevents negative traffic flow when the optimal pricing pattern is replaced by the queueing delay pattern.
  \footnote{
    This condition \eqref{eq:QRP_Condition} easily holds when the slope of the schedule delay cost function for late arrivals is relatively small. 
    However, when the slope is relatively large, this condition may be difficult to satisfy because it requires that bottleneck capacity increase quickly toward the CBD. 
}

\section{Long-term Equilibrium}
\label{sec:LongTermEquilibrium}
Thus far, we have derived the solution to the short-term equilibrium problems.
This section derives the solution to the long-term equilibrium from the solution to the short-term equilibrium.
\par
The long-term equilibrium problem, which determines the distribution of workers regarding the office work ratios and residential locations, can be solved through a constructive approach. 
The key to this approach is a spatial sorting property in the long-term equilibrium state.
The property is that each residential location is classified into the following three types based on the office work ratio of the workers residing there: the office work zone, the remote work zone, and mixed zones.
The office work zone is close to the CBD, where all workers choose the office work ratio of $h=1$. 
Conversely, the remote work zone is situated further away from the CBD, where all workers choose the office work ratio of $h=0$. 
  An equilibrium may exist without any mixed zones depending on parameter settings, but it contains only a single residence if the mixed zone is present.
  The mixed zone is situated between the office work zone and the remote work zone. 
  In this mixed zone, workers choose their office work ratio, $h$, ranging from $0$ to $1$, based on their commuting costs.
  In the remaining portion of this section, we derive the solution to the long-term equilibrium problem and explore this spatial sorting property.
\par
We first rewrite the long-term equilibrium condition using the solution to the short-term equilibrium problem as follows:
\begin{align}
  &\begin{dcases}
    \UE{U}_{i}(h) = \UE{\rho}
  &\mathrm{if}\quad  \UE{Q}_{i}(h) > 0
  \\
  \UE{U}_{i}(h) \leq \UE{\rho} 
  &\mathrm{if}\quad  \UE{Q}_{i}(h) = 0
  \end{dcases}
  &&\forall i \in \ClI^{\ast}, \quad \forall h \in \ClH,
  \\
  &\text{where} \quad  
  \UE{U}_{i}(h)
  = h \left( \theta^{O} - \sum_{j;j\leq i} f_{j} \right)
   - \lambda^{\ast}_{i}(h) 
  + \left( 1 - h \right) \theta^{R} - \UE{r}_{i}.
\end{align}
Because $\lambda^{\ast}_{i}(h)= h\lambda^{\ast}_{i}$ and $\lambda^{\ast}_{i} = \overline{c}(\int_{h \in \ClH}h\UE{Q}_{i}(h), \overline{\mu}_{i})$, we can calculate $\UE{U}_{i}(h)$ as follows:
\begin{align}
  &\UE{U}_{i}(h) = h \left( \theta^{O} - \overline{c}(X_{i}, \overline{\mu}_{i}) - \sum_{j;j\leq i} f_{j} 
      \right) + \left( 1 - h  \right) \theta^{R} - \UE{r}_{i},
\end{align}
  where $X_{i}$ represents the total number of commuters at location $i$, i.e., $X_{i} = \int_{h \in \ClH}h\UE{Q}_{i}(h) \mathrm{d}h$, $0 < X_{i} \leq A_{i}$.
\par
Subsequently, in the equilibrium state, we identify where the mixed zone exists and how many workers living in there choose the office work.
To this end, we introduce a function $G_{i}(X_{i}): \mathbb{R}_{+} \rightarrow \mathbb{R}$ as follows:
\begin{align}
  &G_{i}(X_{i}) \equiv \theta^{O} - \overline{c}(X_{i}, \overline{\mu}_{i}) - \sum_{j;j\leq i} f_{j}
  &&\forall i \in \ClI^{\ast}.
\end{align}
  The value of this function represents the utility gained by workers residing in location $i$ through office work, given that the number of commuters living in $i$ is  $X_{i}$.
From the optimality condition for [ST-SO], we have
  \begin{align}
    &\overline{c}(X_{i}, \overline{\mu}_{i})
    <
    \overline{c}(X_{i+1}, \overline{\mu}_{i+1})
    &&\forall i \in \ClI^{\ast} \setminus \{ i^{\ast} \}, \quad \forall h \in \ClH.
\end{align}
From the inequality, we observe that $G_{i}(X_{i}) > G_{i+1}(X_{i+1})$ for all $i \in \ClI^{\ast} \setminus \{ i^{\ast} \}$ and that the value of $G_{i}(X_{i})$ decreases away from the CBD.
 This observation implies two facts.
 First, as long as $G_{i}(A_{i}) > \theta^{R}$, then workers living closer in the CBD choose the office work ratio $h=1$.
 Second, the mixed zone exists in the smallest index $i^{\ast}$ that satisfies $G_{i}(A_{i}) < \theta^{R}$ because $G_{i}(A_{i})$ decreases away from the CBD.
This is summarized in the following lemma:
\begin{lemma}
  \label{lem:mixde_zone}
    The mixed zone is determined as the smallest index $i^{\ast}$ that satisfies the following inequality:
  \begin{align}
    G_{i}(A_{i}) < \theta^{R}.
  \end{align}
  The residential locations $i<i^{\ast}$ and $i>i^{\ast}$ are the office work zone and the remote work zone, respectively.
\end{lemma}
    Moreover, at the mixed zone $i^{\ast}$, the utility for the office work is equal to that for the remote work because all workers have no incentive to change their office work ratio in the equilibrium state.
    Thus, the number of office workers (commuters) in the mixed zone $X_{i^{\ast}}$ is determined by the following equation to ensure the utility of office work for them is equal to the utility of remote work:\footnote{
        Under a particular parameter situation, \eqref{eq:mixed_zone_G=thetaR} has no non-negative solution for $X_{i^{\ast}}$. 
      In this situation, if commuters exist ($X_{i^{\ast}}>0$), their utilities are lower than those of remote workers.
      Consequently, the mixed zone does not exist because all workers at the location $i^{\ast}$ choose $h=0$. 
      Because this situation is relatively uncommon, this paper focuses on a more general situation where the mixed zone exists.
  }
    \begin{align}
      G_{i^{\ast}}(X_{i^{\ast}}) = \theta^{R}.
      \label{eq:mixed_zone_G=thetaR}
    \end{align}
  From this equation, we can calculate an \textit{aggregate office work ratio} of workers residing in the mixed zone $\eta(i^{\ast}) \equiv X_{i^{\ast}} / A_{i^{\ast}}$ as follows:
\begin{align}
  \eta(i^{\ast}) = \dfrac{1}{A_{i^{\ast}}}  G^{-1}_{i^{\ast}}(\theta^{R}),
  \label{eq:mixed_zone_h}
\end{align}
where $G^{-1}_{i^{\ast}}(c)$ represents the inverse function of $G_{i^{\ast}}(X)$.
  In the equilibrium state, the office work ratio chosen by an \textit{individual} is not uniquely determined because there are innumerable patterns of how individuals choose $h$ such that the aggregate commuting demand is the equilibrium solution $X^{\ast}_{i^{\ast}}$.
  In order to clarify the analysis, without loss of generality, the following analysis treats the equilibrium solution where the office work ratio of individual workers is equal to the aggregate office work ratio, which is one of the equilibrium solutions.\footnote{
     This specification does not affect the welfare analysis in the following section because the commuting costs that affect welfare depend only on the total number of commuters $X_{i^{\ast}}$.  
}
\par
\begin{figure}[tbp]
    \center
       \includegraphics[clip, width=0.50\columnwidth]{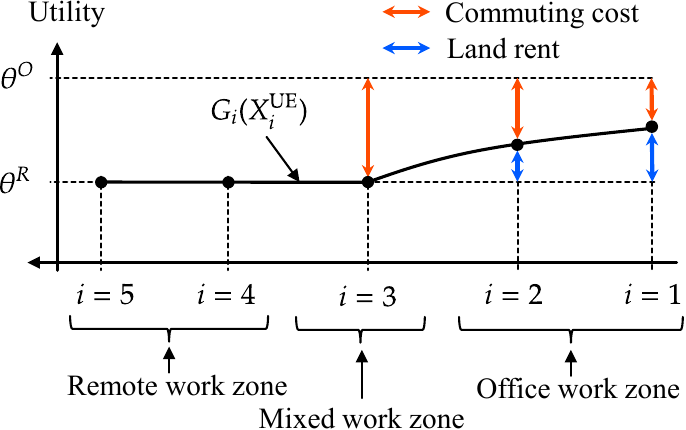}
    \caption{Long-term Equilibrium.}
    \label{fig:LongEquilibriumCost}
  \end{figure}

In the locations further from the mixed zone, i.e., $i > i^{\ast}$, the workers choose $h=0$.
  Summarizing the previous discussion, we can derive the solution to the long-term equilibrium as follows:
\begin{align}
  &\UE{Q}_{i}(h) =
  \begin{dcases}
    A_{i}
    &\mathrm{if}\quad  h = \UE{\eta}(i)
    \\
    0 \quad &\mathrm{otherwise}
  \end{dcases}
  &&\forall i \in \ClI,
  \\
  &\UE{\eta}(i) = 
  \begin{dcases}
    1
    &\mathrm{if}\quad  i < i^{\ast}
    \\
    \dfrac{1}{A_{i^{\ast}}}  G^{-1}_{i^{\ast}}(\theta^{R})
    &\mathrm{if}\quad  i = i^{\ast}
    \\
    0
    &\mathrm{if}\quad  i > i^{\ast}
  \end{dcases}
  \label{eq:eta_i}
\end{align}
  Additionally, we derive the equilibrium utility and land rent. From the long-term equilibrium condition, we have:
\begin{align}
  &\UE{Q}_{i}(h) > 0
  \quad \Rightarrow \quad
  \UE{\rho} = \UE{\eta}(i) \left( \theta^{O} - \overline{c}(\UE{\eta}(i)A_{i}, \overline{\mu}_{i}) - \sum_{j;j\leq i} f_{i} \right)
  + (1-\UE{\eta}(i))\theta^{R}
  \notag
  \\
  &\hspace{120mm} \forall i \in \ClI, \quad \forall h \in \ClH.
  \label{eq:Q_rho=}
\end{align}
For normalization, we set the land rent at the furthest location from the CBD to zero, i.e., $\UE{r}_{I}=0$.
By focusing on \eqref{eq:Q_rho=} for $i=I$, we have $\UE{\rho} = \theta^{R}$.
Based on this equilibrium, we can determine the land rent. Consequently, the solution to the long-term equilibrium problem is as follows:
\begin{proposition}
  The following $\UE{\rho}$, $\{\UE{Q}_{i}(h)\}$, and $\{\UE{r}_{i}\}$ are the solution to the long-term equilibrium problem:
  \begin{align}
    &\UE{\rho} = \theta^{R},
    \\
    &\UE{Q}_{i}(h) =
    \begin{dcases}
      A_{i}
      &\mathrm{if}\quad  h = \UE{\eta}(i)
      \\
      0 \quad &\mathrm{otherwise}
    \end{dcases}
    &&\forall i \in \ClI, \quad \forall h \in \ClH,
    \\
    &\UE{r}_{i} = 
    \begin{dcases}
      G_{i}(\UE{\eta}(i)A_{i}) - \UE{\rho}
      &\mathrm{if}\quad  i < i^{\ast}
      \\
      0
      \quad &\mathrm{otherwise}
    \end{dcases}
    &&\forall i \in \ClI.
  \end{align}
  \label{pro:Solution_LongEquilibrium}
\end{proposition}
\Cref{fig:LongEquilibriumCost} illustrates an example of long-term equilibrium in the case of $I = 5$. 
The horizontal axis represents the distance from the CBD (residential location), while the vertical axis represents utility. 
  The black solid line represents the utility gained by workers residing in location $i$ through office work $G_{i}(\UE{X}_{i})$.
  The red arrows and blue arrows represent the commuting cost and the land rent, respectively.
  In this example, the equilibrium utility is equal to the remote work wage (i.e., $\UE{\rho} = \theta^{R}$), and the mixed zone exists at location $i^{\ast}=3$.
  We observe that the utility gained by office work is equal to the remote work wage in the mixed zone, i.e., $G_{3}(\UE{X}_{3}) = \theta^{R}$.
  Moreover, we find that land rents are positive in the office work zone and zero in the remote work and mixed zones.
 
\section{Welfare Comparison}\label{sec:Comparison}
  In this section, we analyze the impacts of TLC and SWH on the residential location patterns and the commuting patterns of workers. 
  First, we derive the equilibrium in the following four scenarios shown in \Cref{tab:Scenarios} as special cases of the general solution obtained from the integrated equilibrium model described in \Cref{sec:ShortTermEquilibrium,sec:LongTermEquilibrium}. Subsequently, we compare the equilibrium states in these scenarios to analyze the impacts of TLC and SWH. By conducting these comparisons, we can identify the effects of SWH on the long-term equilibrium (residential location pattern), the effects of TLC on the short-term equilibrium (commuting pattern), and the interactions between them.
\begin{table}[tbp]
\center
\caption{Scenarios for the Welfare Analysis.}
\label{tab:Scenarios}
  \begin{tabular}{ |l|c|c| } 
   \hline
   Scenario & $\ClH$ & $K$ \\ 
   \hline 
   (NS):  No schemes & $\ClH=\{1\}$ & $K=1$ \\ 
   (SWH): Staggered work hours & $\ClH=\{1\}$ & $K=2$ \\ 
   (TLC): Telecommuting & $\ClH=[0,1]$ & $K=1$ \\ 
   (CS): Combined scheme & $\ClH=[0,1]$ & $K=2$ \\ 
   \hline
  \end{tabular}
\end{table}
\begin{figure}[tbp]
  \center
    \includegraphics[clip, width=0.50\columnwidth]{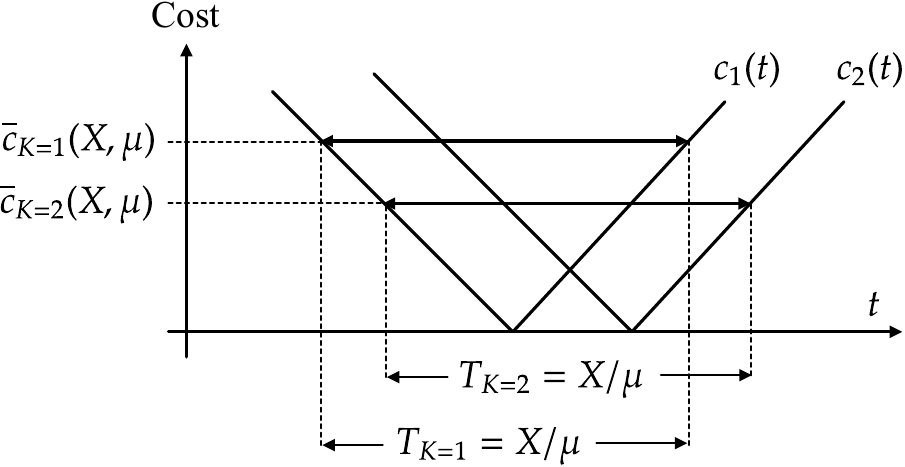}
 \caption{Piecewise Linear Schedule Delay Cost Function.}
 \label{fig:PL_overlinec}
\end{figure}
\par
In the remaining part of this section, we use the notations $\Gamma_{K}(c)$ and $\overline{c}_{K}(X, \mu)$ instead of $\Gamma(c)$ and $\overline{c}(X, \mu)$ because these functions depend on the number of preferred times $K$.
\par
In order to simplify the welfare analysis, we assume a piecewise linear schedule delay cost function, which can be expressed as follows:
\begin{align}
  &c_{k}(t) \equiv
 \begin{dcases}
   \beta (t_{k} - t)
   &\mathrm{if}\quad t < t_{k}
   \\
   \gamma  (t - t_{k})
   &\mathrm{if}\quad t \geq t_{k}
 \end{dcases}
 &&\forall k \in \ClK, \quad \forall t \in \ClT,
\end{align}
where $\beta$ and $\gamma$ are parameters of the schedule delay cost, satisfying $\beta < 1$, and the QRP condition given in ~\eqref{eq:QRP_Condition}.
For this schedule delay cost function $c_{k}(t)$, $\overline{c}_{K}(X, \mu)$ can be calculated as follows:
\begin{align}
  \overline{c}_{K}(X, \mu) = \dfrac{X-\mu d(K-1)}{\mu} \delta
   = \dfrac{X}{\mu} \delta - d(K-1) \delta,
\end{align}
where $\delta = \beta \gamma  / (\beta + \gamma )$ and $d$ is the difference between the preferred arrival times $t_{1}$ and $t_{2}$.
\Cref{fig:PL_overlinec} depicts the relationship between the   
piecewise linear schedule delay cost function $c_{k}(t)$ and its corresponding $\overline{c}_{K}(X, \mu)$.
\par
Let $\TLC{i}{}^{\ast}$ and $\CS{i}{}^{\ast}$ represent the index of the mixed zone in scenarios (TLC) and (CS), respectively, as determined by \Cref{lem:mixde_zone}.
  By adopting a piecewise linear schedule delay cost function, we can derive the following closed-form expression for $\eta(i)$:
\begin{align}
  & \TLC{\eta}(i) 
   = 
  \begin{dcases}
    1 
    &\mathrm{if}\quad  i < \TLC{i}{}^{\ast}
    \\
    \dfrac{\overline{\mu}_{i}}{A_{i}\delta} \left( \theta^{O} - \theta^{R} - \sum_{j;j\leq i} f_{j}\right) 
    &\mathrm{if}\quad  i = \TLC{i}{}^{\ast}
    \\
    0
    &\mathrm{if}\quad  i > \TLC{i}{}^{\ast}
  \end{dcases}
  \\
  &\CS{\eta}(i)
   =
  \begin{dcases}
    1 
    &\mathrm{if}\quad  i < \CS{i}{}^{\ast}
    \\
    \dfrac{\overline{\mu}_{i}}{A_{i}\delta} \left( \theta^{O} - \theta^{R} - \sum_{j;j\leq i} f_{j} + d \delta \right) 
    &\mathrm{if}\quad  i = \CS{i}{}^{\ast}
    \\
    0
    &\mathrm{if}\quad  i > \CS{i}{}^{\ast}
  \end{dcases}
\end{align}
For the sake of simplicity, the following welfare analysis does not focus on the case where all zones are office work zones in the equilibrium states in scenarios (TLC) and (CS).
\footnote{
    For some parameter settings, e.g., when the office work wage $\theta^{O} $ is incredibly high compared to the remote work wage $\theta^{R}$, there can be cases where all zones are office work zones, even in the scenarios (TLC) and (CS).
}
Note that $\NS{\eta}(i) = 1$ and $\SWH{\eta}(i) = 1$ for all $i \in \ClI$, because TLC is not introduced in scenarios (NS) and (SWH).
\par
In the subsequent analysis, we compare a total commuting cost of all workers.
The total commuting cost $TC$ is defined as follows:
  \begin{align}
    TC 
    \equiv 
    \sum_{i \in \ClI} \sum_{k \in \ClK} 
    \int_{h \in \ClH} \int_{t \in \ClT} 
    \lambda_{i}(h) h q_{i,k}(h,t) \mathrm{d} t
    \mathrm{d}h.
    \label{eq:TC}
\end{align}

\begin{figure}[tbp]
  \center
    \includegraphics[clip, width=0.75\columnwidth]{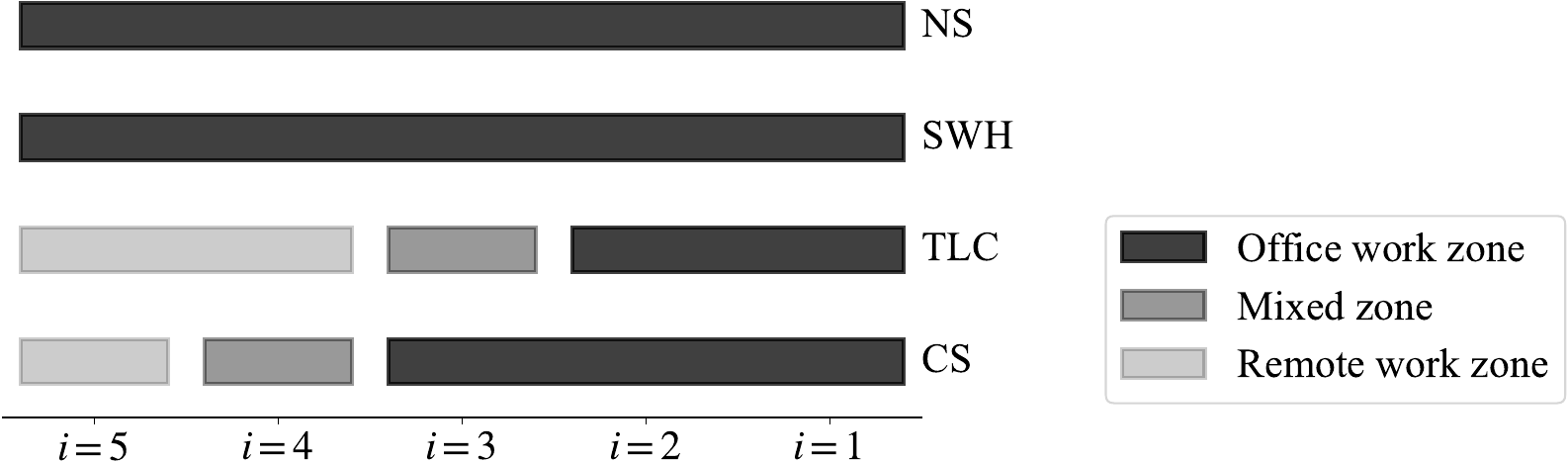}
 \caption{Equilibrium Location Pattern.}
 \label{fig:location_pattern}
\end{figure}

\subsection{Impact of Staggered Work Hours}
This section focuses on the scenarios (NS) and (SWH) to show the impacts of SWH on commuting traffic patterns (short-term effects) and residential location patterns (long-term effects).
\par
  In both scenarios (NS) and (SWH), all workers choose the office work ratio $h=1$. 
  \Cref{fig:location_pattern} depicts the residential pattern in each scenario for the case of $I = 5$ (\Cref{fig:location_pattern} includes the other scenarios for future discussions). 
  The residential pattern in scenario (SWH) is identical to that in scenario (NS). 
  Although all workers commute to the CBD in both scenarios (NS) and (SWH), their commuting costs differ.
  By comparing the equilibrium commuting cost in scenarios (NS) and (SWH), we have the following lemma:
\begin{lemma}
  \label{lem:lambda_NSvsSWH}
    The equilibrium commuting cost in the scenario (SWH) is lower than that in the scenario (NS):
  \begin{align}
    &\SWH{\lambda}_{i}(1) < \NS{\lambda}_{i}(1)
    &&\forall i \in \ClI.
    \label{eq:lambda_NSvsSWH}
  \end{align}
\end{lemma}
\noindent{\bf Proof}.
  See Appendix \ref{subsec:prf_lem_lambda_NSvsSWH}.
\par
\Cref{fig:Cm_NSvsSWH} illustrates the equilibrium commuting cost pattern in scenarios (NS) and (SWH) for the case of $I=5$ as an example.
  At all residential locations, the commuting cost decreases with the introduction of SWH.
  Because this model assumes the piecewise linear schedule cost function and the same choice set of work start times at all locations, the introduction of SWH uniformly reduces commuting costs at each location.
\par
  The introduction of SWH, leading to a uniform decrease in commuting costs, results in the land rent pattern remaining unchanged between scenarios (NS) and (SWH). 
Equilibrium utilities are equal across all locations in both scenarios. 
Because wages (office work wages) are the same at every location in both (NS) and (SWH), the difference in commuting costs between adjacent locations determines the land rent. 
Consequently, the uniform reduction in commuting costs does not alter the land rent pattern. This is encapsulated in the following Lemma.
\begin{lemma}
  \label{lem:r_NSvsSWH}
    The equilibrium land rent in the scenario (SWH) is equal to that in the scenario (NS):
  \begin{align}
    &\SWH{r}_{i} = \NS{r}_{i}
    &&\forall i \in \ClI.
    \label{eq:r_NSvsSWH}
  \end{align}
\end{lemma}
\noindent{\bf Proof}.
  See Appendix \ref{subsec:prf_lem_r_NSvsSWH}.
\par
\Cref{lem:r_NSvsSWH} indicates that the introduction of SWH does not affect the land rent for any location. 
Therefore, the land rent income of the absentee landlord remains unchanged before and after the introduction of SWH.
\par
Furthermore, we have the following theorem:
\begin{theorem}
  \label{thm:SWH_Effect}
    (i) The equilibrium utility for workers in the scenario (SWH) is higher than that in the scenario (NS); 
    (ii) the total commuting cost in the scenario (SWH) is lower than that in the scenario (NS):
  \begin{align}
     \text{(i) } \SWH{\rho} >\NS{\rho},
     \quad
     \text{(ii) } \SWH{TC} < \NS{TC} .
  \end{align}
\end{theorem}
\noindent{\bf Proof}.
  See Appendix \ref{subsec:prf_thm}.

\begin{figure}[tbp]
  \center
  \begin{minipage}[tbp]{0.475\columnwidth}
  \center
    \includegraphics[clip, width=0.95\columnwidth]{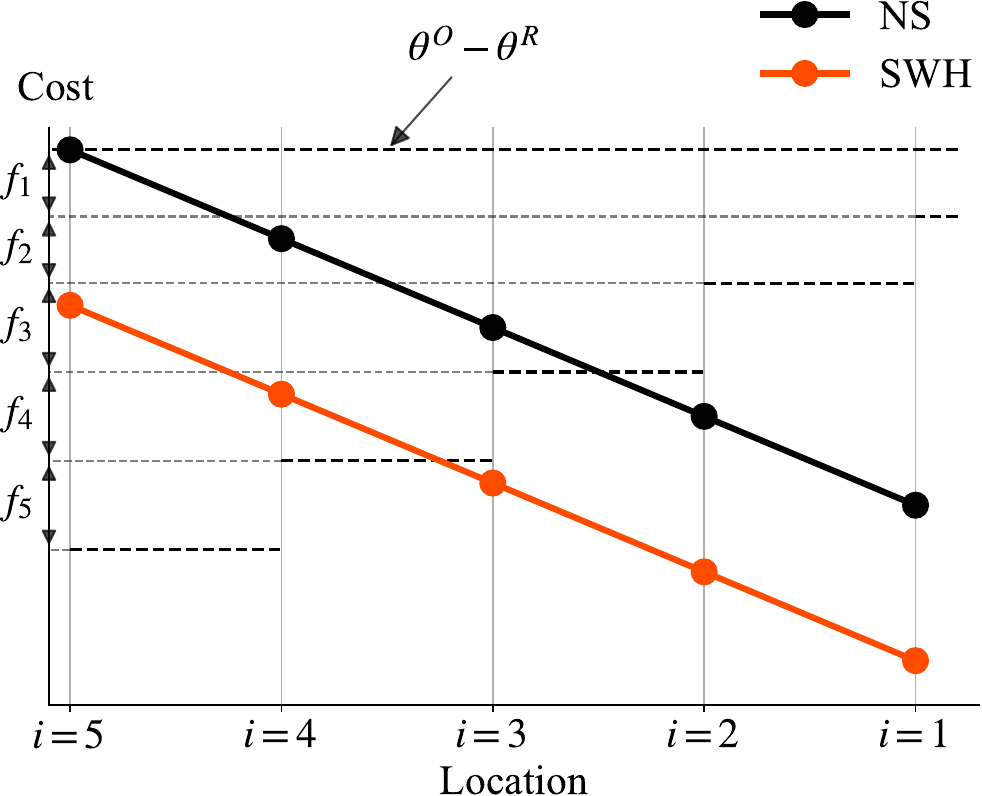}
 \caption{Commuting Cost Pattern; \\ \hspace{15mm} (NS) vs. (SWH).}
 \label{fig:Cm_NSvsSWH}
  \end{minipage}
  \begin{minipage}[tbp]{0.475\columnwidth}
    \center
    \includegraphics[clip, width=0.95\columnwidth]{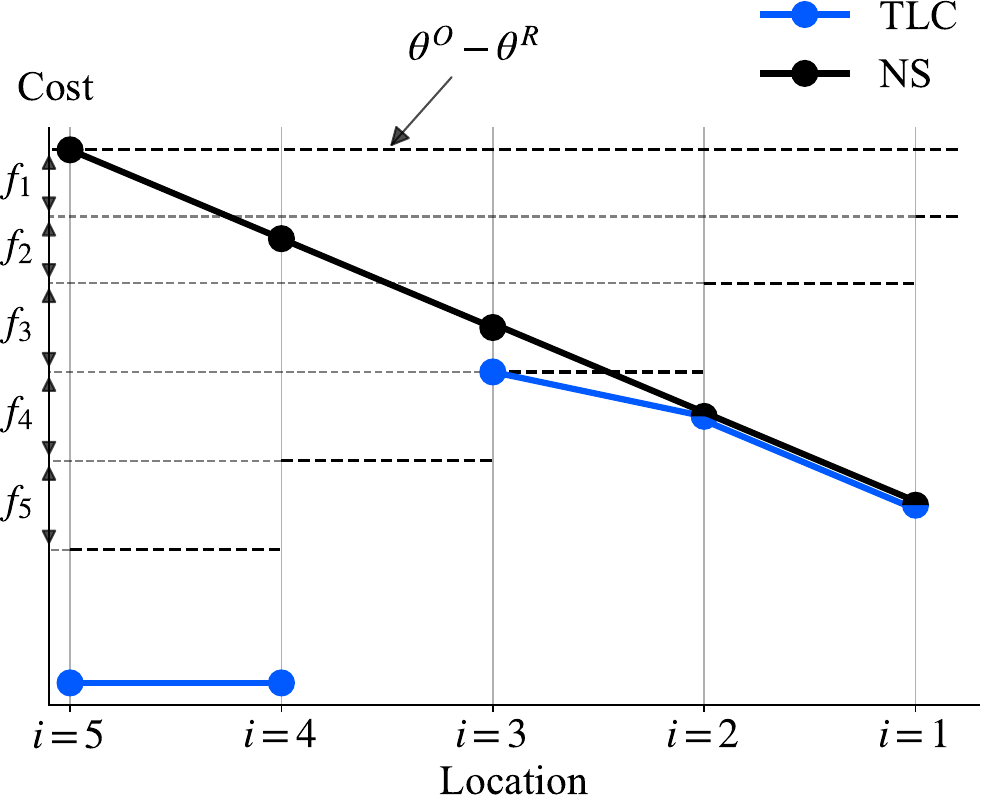}
 \caption{Commuting Cost Pattern; \\ \hspace{15mm} (NS) vs. (TLC).}
 \label{fig:Cm_NSvsTLC}
  \end{minipage}
\end{figure}

\subsection{Impact of Telecommuting}
  In this section, we focus on scenarios (NS) and (TLC) to demonstrate the impacts of TLC on commuting traffic patterns (short-term effects) and residential location patterns (long-term effects). By comparing the equilibrium commuting costs in scenarios (NS) and (TLC), we obtain the following lemma:
\begin{lemma}
  \label{lem:lambda_NSvsTLC}
    For all location $i < \TLC{i}{}^{\ast}$, the equilibrium commuting cost in the scenario (TLC) is equal to that in the scenario (NS):
  \begin{align}
    &\TLC{\lambda}_{i}(\TLC{\eta}(i))= \NS{\lambda}_{i}(1)
    &&\forall i < \TLC{i}{}^{\ast}.
    \label{eq:lambda_NSvsTLC_2}
  \end{align}
    For all location $i \geq \TLC{i}{}^{\ast}$, the equilibrium commuting cost in the scenario (TLC) is lower than that in the scenario (NS):
  \begin{align}
    &\TLC{\lambda}_{i}(\TLC{\eta}(i)) < \NS{\lambda}_{i}(1)
    &&\forall i \geq \TLC{i}{}^{\ast}.
    \label{eq:lambda_NSvsTLC_1}
  \end{align}
\end{lemma}
\noindent{\bf Proof}.
  See Appendix \ref{subsec:prf_lem:lambda_NSvsTLC}.
\par
  \Cref{fig:Cm_NSvsTLC} illustrates the equilibrium commuting cost pattern in scenarios (NS) and (TLC) for the case of $I=5$ as an example.
  In scenario (TLC), the workers residing at location $i<\TLC{i}{}^{\ast}$ (i.e., $i=1, 2$) choose the office work ratio $h=1$, even under TLC.
  Consequently, their commuting costs remain unchanged before and after the introduction of TLC. Conversely, workers who live at location $i > \TLC{i}{}^{\ast}$ (i.e., $i=4, 5$) choose the office work ratio $h = 0$ after the introduction of TLC. This is because the utility from office work, including commuting costs, is lower than that from remote work.
  Workers who reside in the mixed zone adjust their office work ratio to ensure their utility is not lower than that from remote work.
\par
\Cref{fig:location_pattern} shows the residential location pattern in scenarios (NS) and (TLC) for the case of $I=5$ as an example.
The residential pattern is classified into three types of locations: the mixed zone, the office worker zone, and the remote worker zone.
  Consequently, the land rent pattern varies with the introduction of TLC, as stated in the following lemma:
\begin{lemma}
  \label{lem:r_NSvsTLC}
    The equilibrium land rent in the scenario (TLC) is equal to or lower than that in the scenario (NS):
  \begin{align}
    &\TLC{r}_{i} \leq \NS{r}_{i} 
    &&\forall i \in \ClI.
    \label{eq:r_NSvsTLC} 
  \end{align}
\end{lemma}
\noindent{\bf Proof}.
  See Appendix \ref{subsec:prf_lem:r_NSvsTLC}.
\par
  \Cref{lem:r_NSvsTLC} indicates that the introduction of TLC leads to a decrease in land rent for all locations. This implies that the land rent income of the absentee landlord decreases after the introduction of TLC.
\par
Furthermore, we have the following theorem:
\begin{theorem}
  \label{thm:RemoteWorkEffect}
    (i) The equilibrium utility for workers in the scenario (TLC) is higher than that in the scenario (NS);
    (ii) the total commuting cost in the scenario (TLC) is lower than that in the scenario (NS):
  \begin{align}
    \text{(i) } \TLC{\rho} > \NS{\rho}, 
    \quad 
    \text{(ii) } \TLC{TC} < \NS{TC}.
  \end{align}
\end{theorem}
\noindent{\bf Proof}.
  See Appendix \ref{subsec:prf_thm}.
\par
The effects of introducing TLC is similar to those of introducing SWH (\Cref{thm:SWH_Effect}).
However, the land rent decreases with the introduction of TLC, whereas it does not change with the introduction of SWH.

\begin{figure}[tbp]
  \center
  \begin{minipage}[tbp]{0.475\columnwidth}
    \center
    \includegraphics[clip, width=0.95\columnwidth]{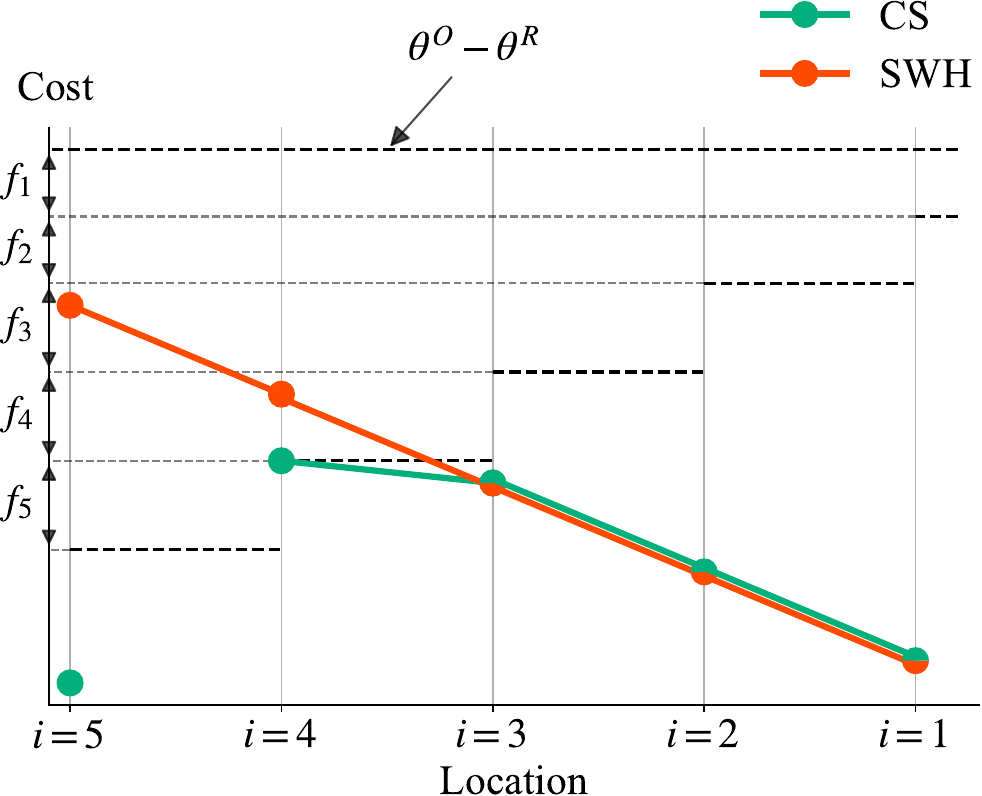}
    \caption{Commuting Cost Pattern; \\ \hspace{18mm} (SWH) vs. (CS).}
    \label{fig:Cm_SWHvsCS}
  \end{minipage}
  \begin{minipage}[tbp]{0.475\columnwidth}
  \center
    \includegraphics[clip, width=0.95\columnwidth]{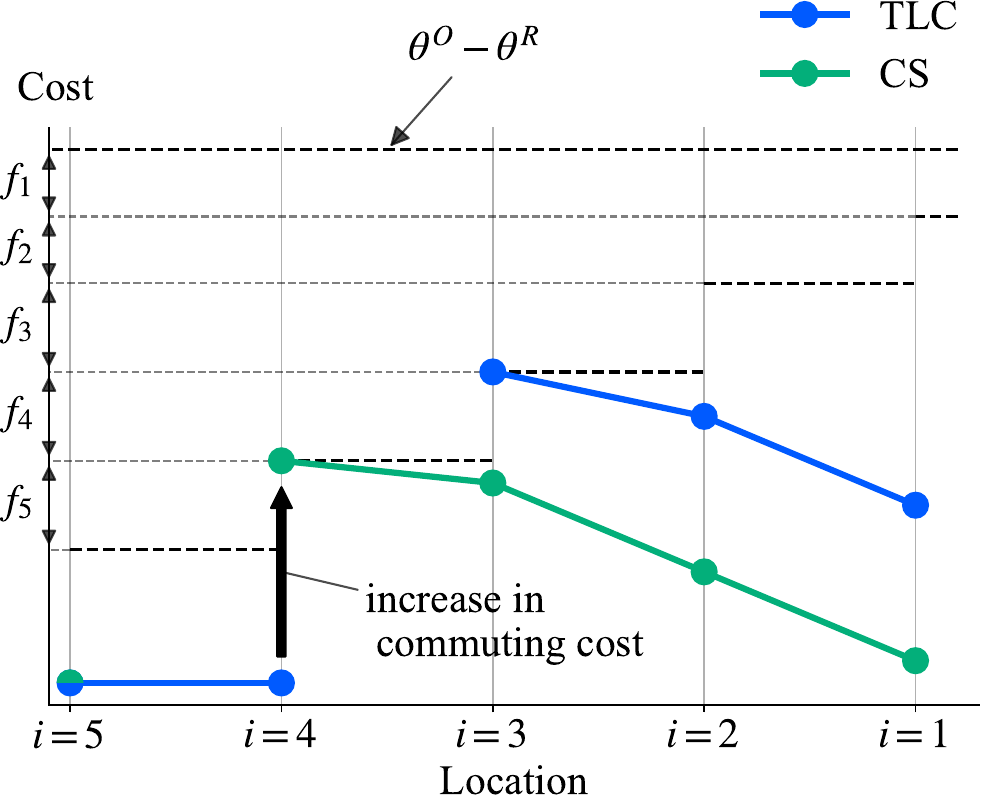}
 \caption{Commuting Cost Pattern; \\ \hspace{18mm} (TLC) vs. (CS).}
 \label{fig:CmmCst_TLCvsCS}
  \end{minipage}
\end{figure}

\subsection{Impacts of Combined Scheme}
We investigate the impacts of TLC and SWH using two comparisons: i) scenario (SWH) vs (CS) and ii) scenario (TLC) vs (CS).
These comparisons enable us to find the impacts of introducing TLC in addition to SWH and introducing SWH in addition to TLC, respectively.

\subsubsection{Combined Scheme vs. Staggered Work Hours.}
By comparing the equilibrium commuting cost in scenarios (SWH) and (CS), we have the following lemmas regarding the equilibrium commuting cost and land rent:
\begin{lemma}
  \label{lem:lambda_SWHvsCS}
    For all locations $i < \CS{i}{}^{\ast}$, the equilibrium commuting cost in the scenario (CS) is equal to that in the scenario (SWH):
  \begin{align}
    &\CS{\lambda}_{i}(\CS{\eta}(i)) = \SWH{\lambda}_{i}(1)
    &&\forall i < \CS{i}{}^{\ast}.
    \label{eq:lambda_SWHvsCS_2}
  \end{align}
    For all locations $i \geq \CS{i}{}^{\ast}$, the equilibrium commuting cost in the scenario (CS) is lower than that in the scenario (SWH):
  \begin{align}
    &\CS{\lambda}_{i}(\CS{\eta}(i)) < \SWH{\lambda}_{i}(1)
    &&\forall i \geq \CS{i}{}^{\ast}.
    \label{eq:lambda_SWHvsCS_1}
  \end{align}
\end{lemma}
\noindent{\bf Proof}.
  See Appendix \ref{subsec:prf_lem:lambda_SWHvsCS}.
\begin{lemma}
  \label{lem:r_SWHvsCS}
    The equilibrium land rent in the scenario (CS) is equal to or lower than that in the scenario (SWH):
  \begin{align}
    &\CS{r}_{i} \leq \SWH{r}_{i}
    &&\forall i \in \ClI.
    \label{eq:r_SWHvsCS}
  \end{align}
\end{lemma}
\noindent{\bf Proof}.
  See Appendix \ref{subsec:prf_lem:r_SWHvsCS}.
\par
\Cref{fig:Cm_SWHvsCS} and \Cref{fig:location_pattern} depict the equilibrium commuting cost and location patterns in scenarios (SWH) and (CS) for the case of $I=5$ as an example.
As shown in \Cref{fig:location_pattern}, the mixed zone is positioned further from the CBD in scenario (CS) compared to scenario (SWH).
Thus, the commuting costs do not change at the locations $i=1,2,3$ and 
decrease at the locations $i=4,5$, with the introduction of TLC in addition to SWH.
\par
By \Cref{lem:lambda_SWHvsCS,lem:r_SWHvsCS}, we have the following theorem:
\begin{theorem}
  \label{thm:SWHvsCS}
    (i) The equilibrium utility for workers in the scenario (CS) is higher than that in the scenario (SWH);
    (ii) the total commuting cost in the scenario (CS) is lower than that in the scenario (SWH):
  \begin{align}
    \text{(i) }\CS{\rho} > \SWH{\rho},
    \quad 
    \text{(ii) }\CS{TC} < \SWH{TC}.
  \end{align}
\end{theorem}
\noindent{\bf Proof}.
  See Appendix \ref{subsec:prf_thm}.
\par
\Cref{thm:SWHvsCS} suggests that the introduction of TLC alongside SWH results in lower total traffic costs and higher utility for workers than implementing SWH alone.

\begin{table}[tbp]
  \center
  \caption{Summary of \Cref{lem:lambda_TLCvsCS,lem:r_TLCvsCS}.}
  \label{tab:summary}
    \begin{tabular}{|c|c|c|c|c|} 
     \hline
     Zone type in (TLC) & 
     Zone type in (CS) & 
     Commuting cost  & 
     Land rent  &
     Wage 
     \\
     {} & 
     {} & 
     $\lambda_{i}(h)$ & 
     $r_{i}$ &
     $h\theta^{O} + (1-h) \theta^{R}$
     \\ 
     \hline \hline 
     Office work zone & Office work zone
     & $\searrow$ & $\nearrow$ & $\longrightarrow$
     \\
     Mixed zone& Office work zone 
     & $\searrow$ & $\nearrow$ & $\nearrow$ 
     \\
     Remote work zone& Office work zone
     & $\nearrow$ & $\nearrow$ & $\nearrow$ 
     \\
     Remote work zone& Mixed zone
     & $\nearrow$ & $\longrightarrow$ & $\nearrow$ 
     \\ 
     Remote work zone& Remote work zone
     & $\longrightarrow$ & $\longrightarrow$ & $\longrightarrow$
     \\
     \hline
    \end{tabular}
  \end{table}

\subsubsection{Combined Scheme vs. Telecommuting.}
  We compare scenarios (TLC) and (CS) and show a paradoxical fact. 
  We first focus on the equilibrium residential location pattern.
  As shown in \Cref{fig:location_pattern}, the office work zone expands with the introduction of SWH in addition to TLC,i.e., $\CS{i}{}^{\ast} > \TLC{i}{}^{\ast}$. The expansion of office work zones is attributed to the decrease in commuting costs associated with the introduction of SWH. Consequently, locations where workers are willing to engage in office work extend further away from the CBD. 
  From these facts, we have the following lemmas regarding the equilibrium commuting cost and land rent:
\begin{lemma}
  \label{lem:lambda_TLCvsCS}
    For all locations $i \leq \TLC{i}{}^{\ast}$, the equilibrium commuting cost in the scenario (CS) is lower than that in the scenario (TLC):
  \begin{align}
    &\CS{\lambda}_{i}(\CS{\eta}(i)) \leq \TLC{\lambda}_{i}(\TLC{\eta}(i)) 
    &&\forall i < \TLC{i}{}^{\ast}.
    \label{eq:lambda_TLCvsCS_1}
  \end{align}
    For all locations $i > \TLC{i}{}^{\ast}$, the equilibrium commuting cost in the scenario (CS) is higher than that in the scenario (TLC):
  \begin{align}
    &\CS{\lambda}_{i}(\CS{\eta}(i)) > \TLC{\lambda}_{i}(\TLC{\eta}(i))
    &&\forall i > \TLC{i}{}^{\ast}.
    \label{eq:lambda_TLCvsCS_2}
  \end{align}
\end{lemma}
\noindent{\bf Proof}.
  See Appendix \ref{subsec:prf_lem:lambda_TLCvsCS}.
\begin{lemma}
  \label{lem:r_TLCvsCS}
    For all locations $ i < \CS{i}{}^{\ast}$, the equilibrium land rent in the scenario (CS) is higher than that in the scenario (TLC):
  \begin{align}
    &\CS{r}_{i} > \TLC{r}_{i}
    &&\forall i < \CS{i}{}^{\ast}.
  \end{align}
  For all locations $i \geq \CS{i}{}^{\ast}$, the equilibrium land rent in the scenario (CS) is equal to that in the scenario (TLC):
  \begin{align}
    &\CS{r}_{i} = \TLC{r}_{i}
    &&\forall i \geq \CS{i}{}^{\ast}.
  \end{align}
\end{lemma}
\noindent{\bf Proof}.
  See Appendix \ref{subsec:prf_lem:r_TLCvsCS}.
\par
  \Cref{fig:CmmCst_TLCvsCS} and \Cref{fig:location_pattern} depict the equilibrium commuting cost and location patterns in scenarios (TLC) and (CS) for the case of $I = 5$ as an example.
  By introducing SWH in addition to TLC, the commuting cost decreases at the locations $i=1,2,3$; however, the commuting cost increases at the location $i=4$ due to the rise in the number of office workers (commuting demand).
\par
\Cref{tab:summary} summarizes the statements of \Cref{lem:lambda_TLCvsCS,lem:r_TLCvsCS}.
There are both locations where commuting costs increase and decrease.
Therefore, it is not possible to assess in general whether scenario (TLC) or (CS) has a higher total commuting cost, unlike the other comparisons of the scenarios.
Because the office work zone is larger in Scenario (CS) than in Scenario (TLC), the total commuting cost is higher in Scenario (CS), depending on the parameters.
This is summarized in the following theorem:
\begin{theorem}
  \label{thm:RemoteWorkParadox}
    (i) The equilibrium utility for workers in the scenario (CS) is equal to that in the scenario (TLC);
    (ii) the total commuting cost in the scenario (CS) may be higher than that in the scenario (TLC):
  \begin{align}
    \text{(i) }\CS{\rho} = \TLC{\rho}.
  \end{align}
  \end{theorem}
  \noindent{\bf Proof}.
  See Appendix \ref{subsec:prf_thm} and \ref{sec:ExampleoftheParadox}.
  Appendix \ref{sec:ExampleoftheParadox} provides a specific example of cases where the total commuting cost in the scenario (CS) is higher than that in the scenario (TLC).
  \par
  \Cref{thm:RemoteWorkParadox} implies a paradoxical phenomenon in which introducing SWH in addition to TLC may yield higher total commuting costs than only introducing TLC despite the lack of improvement in worker utility.
  This paradox arises due to induced demand caused by reduced commuting costs due to SWH. 
  Although SWH spreads out preferred arrival times, additional traffic congestion arises around each preferred arrival time because the equilibrium state balances the utility of office workers with that of remote workers. 
  Thus, the paradox arises when the increase in total commuting costs due to this additional demand offsets the decrease in total commuting costs of existing office workers due to the introduction of SWH.
  \par
    This paradox is more likely to occur when there are few office workers in the scenario (TLC),  i.e., the office work zone is small in the scenario (TLC).
    This is because only some workers directly benefit from SWH in such a situation.
    The size of the office work zone is smaller as the difference between remote work wages and office work wages $(\theta^{O} - \theta^{R})$ is smaller (\Cref{lem:mixde_zone}).
    Thus, when the wage of remote work is as high as that of office work, this paradox is more likely to occur.

\section{Conclusions}
  In this study, we investigated the long- and short-term effects of TLC, SWH, and their combined scheme on peak-period congestion and location patterns. 
  In order to facilitate a unified comparison of the long- and short-term effects of these schemes, we developed a new equilibrium analysis approach that integrates the long-term equilibrium (location and percentage of telecommuting choice) and short-term equilibrium (preferred and actual arrival times choice) while considering their interaction. 
  We derived closed-form solutions for the long- and short-term equilibrium using their mathematical structures, the same as optimal transport problems. As expected, these solutions reveal the positive impact of introducing TLC or SWH on worker utility. However, to our surprise, we discovered a paradox: introducing SWH alongside TLC does not improve workers' utility compared to a scenario with TLC alone. 
  Moreover, we showed that this paradox could lead to higher total commuting costs.
  \par
    Several key issues could be addressed to advance this study's findings further. 
    First, we should work on a comprehensive sensitivity analysis of parameters such as bottleneck capacity pattern, land area pattern, and wages. 
    We also need to investigate how the results depend on the assumptions in this study to derive the closed-form solution, including homogeneous workers, bottleneck capacity pattern, and the schedule delay cost function. 
    These analyses would clarify quantitatively and in more detail when the paradox is more likely to occur.
    As more applied future studies, we should analyze the long-term and short-term effects of TLC and other transportation demand management (TDM) policies other than SWH. 
    This study revealed that reducing individual commuting costs through SWH may lead to induced demand and a paradox when TLC is in place. 
    Hence, TDM policies other than SWH might also result in similar paradoxes. 
    It is also significant to examine the relevance and similarity of other paradoxes in the bottleneck model, such as capacity expansion paradoxes \citep{Arnott1993-mq,He2022-ca}. 
    Addressing these issues will enhance our understanding of the intricate interactions between multiple policies.
  
\section*{Acknowledgement}
  The authors are grateful to the associate editor and three anonymous referees for valuable comments on an earlier version of the paper. 
  This work was supported by Council for Science, Technology and Innovation (CSTI), Cross-ministerial Strategic Innovation Promotion Program (SIP), the 3rd period of SIP ``Smart Infrastructure Management System'' Grant Number JPJ012187 (Funding agency: Public Works Research Institute).
  This work was supported by JSPS KAKENHI Grant Numbers JP20J21744, JP21H01448, JP24K00999, JP20K14843, and JP23K13418.

\appendix
\section{List of Notation}
\label{sec:List_Notation}
\Cref{tab:List_Notation_1}, \Cref{tab:List_Notation_2}, and \Cref{tab:List_Notation_3} list notation frequently used in this paper.
\begin{table}[!ht]
  \centering
  \caption{List of Notation (1).}
  \label{tab:List_Notation_1}
  \begin{tabularx}{\linewidth}{lX}
    \textit{Parameters}  & {}
    \\
    \hline \hline
      $N$ & Number of work days in a unit term
       \\
       $I$ & Number of origin nodes/bottlenecks in the corridor network
       \\
       $\ClI = \{ i \mid i=1,2,...,I \}$
       & Set of origin nodes/bottlenecks
       \\
       $A_{i}$ & Area of residential location $i$
       \\
       $f_{i}$ & Free flow travel time from origin $i$ to the destination
       \\
       $\mu_{i}$ & Capacity of bottleneck $i$
       \\ 
       $\overline{\mu}_{i} = \mu_{i} - \mu_{(i+1)}$ &
       Difference between the capacities of bottleneck i and the upstream bottleneck $i + 1$, where $\mu_{(N+1)} = 0$
       \\
       $\ClT$ & Set of arrival times (assignment period)
       \\
       $K$ & Number of official work start times (preferred arrival times)
       \\
       $\ClK = \{ k \mid k=1,2,...,K \}$ & Set of in index of preferred arrival times
       \\
       $t_{k}$ & $k$ th preferred arrival time 
       \\
       $d_{k} = t_{k+1} - t_{k}$ & Difference between preferred arrival times 
       \\
       $c(T)$ & Base schedule delay cost function
       \\
       $c_{k}(t)$ & Schedule delay cost function of the $k$ th preferred arrival time
       \\
       $\alpha$ &
        Value of time (VOT) of all commuters ($\alpha = 1$)
       \\
       $\widehat{c}(t)$, $\widehat{\ClT}_{k}$, $\Gamma(c)$ and $\overline{c}(X, \mu)$ &
       Functions and sets related to $c_{k}(t)$ defined by \eqref{eq:widehat_c}, \eqref{eq:widehat_T} \eqref{eq:gamma_c}, and \eqref{eq:overline_c}, respectively
       \\
       $\theta^{O}$, $\theta^{R}$&
       Wage of office work and remote work, respectively ($\theta^{O} > \theta^{R}$)
       \\
       $h$ & Office work ratio
       \\
       $\ClH=[0,1]$ & Set of office work ratios
       \\
       $\beta$, $\gamma$ &
        Parameters of the schedule delay cost function, satisfying $\beta < 1$.
       \\
       \hline
       \\
      \end{tabularx}%
    \end{table}
      
   \begin{table}[!ht]
  \centering
  \caption{List of Notation (2).}
  \label{tab:List_Notation_2}
      \begin{tabularx}{\linewidth}{lX}
        \textit{Variables} & {}
        \\
        \hline \hline
       $\tau_{i}(t)$ &
       Arrival time at bottleneck $i$ for commuters whose destination arrival time is $t$
       \\
       $\sigma_{i}(t)$ &
        Departure time at bottleneck $i$ for commuters whose destination arrival time is $t$
       \\
       $q_{i,k}(h, t)$ &
        Arrival flow rate at the CBD of $(h,i,k,t)$-worker
       \\
       $w_{i}(t)$ & Queuing delay at the bottleneck $i$ for the commuters whose destination arrival time is $t$
       \\
       $p_{i}(t)$ & Optimal price at the bottleneck $i$ for the commuters whose destination arrival time is $t$
       \\
       $Q_{i}(h)$ & Number of workers who choose the office work ratio $h$ and location $i$
       \\
       $r_{i}$ & Land rent at the location $i$
       \\
       $\lambda_{i}(h)$ & Equilibrium commuting cost of $(h,i)$-workers
       \\
       $\rho_{i,k}$ & Equilibrium utility of $(i,k)$-workers
       \\
       $X_{i}$ & Commuting demand in the location $i$, i.e., $X_{i} \equiv \int_{h \in \ClH} h Q_{i}(h) \mathrm{d}h$
       \\
       $\eta(i)$ & Aggregate office work ratio at the location $i$
       \\
       \hline
       \\
      \end{tabularx}%
   \end{table}         
   \begin{table}[!ht]
  \centering
  \caption{List of Notation (3).}
  \label{tab:List_Notation_3}
    \begin{tabularx}{\linewidth}{lX}
      \textit{Superscripts meaning of variable $X$} & {}
      \\
      \hline \hline
        $\SO{X}$ & $X$ in the dynamic system optimal state (DSO)
        \\
        $\UE{X}$ & $X$ in the dynamic user equilibrium state (DUE)
        \\
        $\NS{X}$ & $X$ in the equilibrium under the scenario (NS)
        \\
        $\SWH{X}$ & $X$ in the equilibrium under the scenario (SWH)
        \\
        $\TLC{X}$ & $X$ in the equilibrium under the scenario (TLC)
        \\
        $\CS{X}$ &  $X$ in the equilibrium under the scenario (CS)
        \\
        \hline
     \end{tabularx}%
   \end{table}%

\newpage
\section{Proof of Theorems and Lemmas}
We prove the theorems and lemmas in this section.
\subsection{Proof of \Cref{lem:q=mu}}\label{subsec:proof_lem_q=mu}
  From \Cref{asm:p>0q>0}, we find  
  \begin{align}
    &\sum_{k \in \ClK} \int_{h \in \ClH}
      h \SO{q}_{h,I,k}(t) \mathrm{d}h= \mu_{I}
    &&\forall t \in \SO{\ClT}_{I}.
  \end{align}
  Moreover, in the case of $i=I-1$, the following relationship holds:
  \begin{align}
    \sum_{k \in \ClK}
    \int_{h \in \ClH}h \SO{q}_{h,I-1,k}(t)\mathrm{d}h = \mu_{I-1} - \mu_{I}
    \quad \forall t \in \SO{\ClT}_{I-1}.
  \end{align}
  Thus, by applying this recursively, the bottleneck capacity constraint can be rewritten as follows:
  \begin{align}
    &\sum_{k \in \ClK} \int_{h \in \ClH} 
    q_{i,k}(h,t) \mathrm{d}h \leq \mu_{i} - \mu_{i+1}
    &&\forall i \in \ClI^{\ast}, \  \forall t \in \ClT.
  \end{align}
This completes the proof. \hfill \Halmos

\subsection{Proof of \Cref{lem:Solution_BottleneckBasedProblem}}
\label{subsec:proof_lem_Solution_BottleneckBasedProblem}
This subsection proves \Cref{lem:Solution_BottleneckBasedProblem} by deriving the solution to the single bottleneck problem.
Owing to this, we introduce the short-term optimal problem in the single bottleneck case as follows:
\begin{align}
  &\text{[Single]} \notag
  \\
  &\min_{\{ q_{k}(h,t) \} \geq 0}. \quad 
  \sum_{k \in \ClK} \int_{h \in \ClH}\int_{t \in \ClT} h  c_{k}(t) q_{k}(h,t) \mathrm{d}t \mathrm{d}h
  \\
  &\text{s.t.} \quad 
  \sum_{k \in \ClK} \int_{t \in \ClT} q_{k}(h,t) 
  \mathrm{d} t = Q(h)
  &&\forall h \in \ClH
  \quad [\lambda(h)],
  \\
  &\quad \quad \sum_{k \in \ClK} \int_{h \in \ClH} 
  h  q_{k}(h,t) \mathrm{d}h \leq \mu
  &&\forall t \in \ClT
  \quad [w(t)].
  \label{eq:Single_EOP_BNcapa}
\end{align}
The variables in [], which represent the Lagrange multipliers corresponding to those constraints.
We show the optimality condition for [Single] as follows:
  \begin{align}
    &\sum_{k \in \ClK} \int_{t \in \ClT} q^{\ast}_{k}(h,t) 
      \mathrm{d} t = Q(h)
    &&\forall h \in \ClH,
    \\
    &\begin{dcases}
      h  c_{k}(t) + h  w^{\ast}(t)
      = \lambda^{\ast}(h)
      &\mathrm{if}\quad  q^{\ast}_{k}(h,t) > 0
      \\
      h  c_{k}(t) + h  w^{\ast}(t)
      \geq \lambda^{\ast}(h)
      &\mathrm{if}\quad  q^{\ast}_{k}(h,t) = 0
    \end{dcases}
    &&\forall k \in \ClK, \quad \forall t \in \ClT, \quad \forall h \in \ClH,
    \\
    &\begin{dcases}
      \sum_{k \in \ClK} \int_{h \in \ClH} 
      h  q^{\ast}_{k}(h,t) \mathrm{d}h = \mu 
      &\mathrm{if}\quad  w^{\ast}(t) > 0
      \\
      \sum_{k \in \ClK} \int_{h \in \ClH} 
      h  q^{\ast}_{k}(h,t) \mathrm{d}h \leq \mu 
      &\mathrm{if}\quad  w^{\ast}(t) = 0
    \end{dcases}
    &&\forall t \in \ClT.
  \end{align}
  In order to clarify the mathematical structure of [Single], we replace the decision variables of the problem from $q_{i,k}(h,t)$ to $\widehat{q}_{h,i,k}(t)$, where $\widehat{q}_{h,i,k}(t) \equiv h  q_{i,k}(h,t)$.
Using this, the problem [Single] is transformed into the following problem:
\begin{align}
  &\text{[Single']} \notag
  \\
  &\min_{\{\widehat{q}_{k}(h,t)\} \geq 0}. \quad 
  \sum_{k \in \ClK} \int_{h \in \ClH}  \int_{t \in \ClT} 
  c_{k}(t) \widehat{q}_{k}(h,t) \mathrm{d}t\mathrm{d}h
   \\
   &\text{s.t.} \quad 
   \sum_{k \in \ClK} \int_{t \in \ClT} \widehat{q}_{k}(h,t) 
   \mathrm{d} t = h Q(h)
   &&\forall h \in \ClH,
   \\
   &\quad \quad \sum_{k \in \ClK} \int_{h \in \ClH}  
   \widehat{q}_{k}(h,t) \mathrm{d}h \leq \mu
   &&\forall t \in \ClT.
\end{align}
The optimality condition for [Single'] is as follows:
\begin{align}
  &\sum_{k \in \ClK} \int_{t \in \ClT} \widehat{q}^{\ast}_{k}(h,t) 
  \mathrm{d} t = h Q(h)
  &&\forall h \in \ClH,
\\
&\begin{dcases}
  c_{k}(t) + \widehat{p}^{\ast}(t)
  = \widehat{\lambda}^{\ast}(h)
  &\mathrm{if}\quad  \widehat{q}^{\ast}_{k}(h,t) > 0
  \\
  c_{k}(t) + \widehat{p}^{\ast}(t)
  \geq \widehat{\lambda}^{\ast}(h)
  &\mathrm{if}\quad  \widehat{q}^{\ast}_{k}(h,t) = 0
\end{dcases}
&&\forall k \in \ClK, \quad \forall t \in \ClT, \quad \forall h \in \ClH,
\\
&\begin{dcases}
  \int_{h \in \ClH} \sum_{k \in \ClK} 
  \widehat{q}^{\ast}_{k}(h,t) = \mu 
  &\mathrm{if}\quad  \widehat{p}^{\ast}(t) > 0
  \\
  \int_{h \in \ClH} \sum_{k \in \ClK} 
  \widehat{q}^{\ast}_{k}(h,t) \leq \mu 
  &\mathrm{if}\quad  \widehat{p}^{\ast}(t) = 0
\end{dcases}
&&\forall t \in \ClT.
\end{align}
The objective function of the problem [Single'] is independent of the office work ratio $h \in \ClH$.
Because the problems [Single'] and [Single] are equivalent, the equilibrium commuting cost is determined as $\widehat{\lambda}^{\ast}(h) = \widehat{\lambda}^{\ast}$ where $\widehat{\lambda}^{\ast}$ is calculated as follows:
\begin{align}
  \widehat{\lambda}^{\ast} 
  = \overline{c}\left(\int_{h \in \ClH} h Q(h) \mathrm{d}h, \mu\right),
\end{align}
where $\overline{c}(X, \mu)$ is defined in \eqref{eq:overline_c}.
The function $\overline{c}(X, \mu)$ is convex of $X$.
Using these functions, we obtain the solution to [Single'] as follows:
  \begin{align}
    &\widehat{\lambda}^{\ast}(h) 
    = \widehat{\lambda}^{\ast} 
    &&\forall h \in \ClH, 
    \\
    &\widehat{p}^{\ast}(t) = \max \left\{ 0, \widehat{\lambda}^{\ast} - \widehat{c}(t)\right\}
    &&\forall t \in \ClT,
    \\
    &\int_{h \in \ClH} \widehat{q}^{\ast}_{k}(h,t) = 
    \begin{dcases}
      \mu
      &\mathrm{if}\quad  t \in \Gamma(\widehat{\lambda}^{\ast}(h)) \cap \widehat{\ClT}_{k}
      \\
      0
      &\mathrm{otherwise}
    \end{dcases}
    &&\forall k \in \ClK, \quad \forall t \in \ClT.
  \end{align}
From the solution to [Single'], we have the solution to [Single] as follows.
The variable set $\{ q^{\ast}_{k}(h, t) \}$, $\{ \lambda^{\ast}(t) \}$ and $\{ p^{\ast}(t) \}$ is the solution to the problem [Single].
\begin{align}
  &\lambda^{\ast}(h) = h  \lambda^{\ast}
  &&\forall h \in \ClH,
  \\
  &p^{\ast}(t) = \max \left\{ 0, \lambda^{\ast} - \widehat{c}(t)\right\}
  &&\forall t \in \ClT,
  \\
  &\int_{h \in \ClH} q^{\ast}_{k}(h,t) \mathrm{d}h= 
  \begin{dcases}
    \mu
    &\mathrm{if}\quad  t \in \Gamma(\lambda^{\ast}) \cap \widehat{\ClT}_{k}
    \\
    0
    \quad &\mathrm{otherwise}
  \end{dcases}
  &&\forall k \in \ClK, \quad \forall t \in \ClT.
\end{align}
This completes the proof of \Cref{lem:Solution_BottleneckBasedProblem}.
\hfill \Halmos

\subsection{Proof of \Cref{pro:QRP,pro:DUE_Solution}}
\label{subsec:prf_QRP_DUE_Solution}
  The sufficient condition that variables shown in \Cref{pro:DUE_Solution} are the solution to the short-term equilibrium problem is that the variables satisfy the departure time choice condition \eqref{eq:DUE_DeparturePrefferd_Time_Choice}, the queueing condition \eqref{eq:DUE_QueueingCondition}, and the demand conservation condition \eqref{eq:DUE_CommuterCnsv}.
  From \eqref{eq:DUE_solution_q_sum}, we can immediately confirm that the demand conservation condition \eqref{eq:DUE_CommuterCnsv} is satisfied.
  Moreover, because the variables $\{ \SO{p}_{i}(t) \}$ and $\{ \SO{\lambda}_{i}(h) \}$ are the solution to the short-term optimal problem, they satisfy the optimality condition for the short-term optimal problem.
  Thus, the departure time choice condition \eqref{eq:DUE_DeparturePrefferd_Time_Choice} are satisfied.
  Therefore, if the variables satisfy the queueing condition \eqref{eq:DUE_QueueingCondition} and the demand conservation condition \eqref{eq:DUE_CommuterCnsv}, the variables are the solution to the short-term equilibrium problem.
  \par
  We first check the queueing condition \eqref{eq:DUE_QueueingCondition}.
  By summing \eqref{eq:DUE_solution_q} recursively, we have 
  the aggregate arrival rate as follows:
  \begin{align}
    &\sum_{j;j\geq i} 
    \sum_{k \in \ClK}
    \int_{h \in \ClH} h \UE{q}_{j,k}(h,t) \mathrm{d}h
    = 
    \begin{dcases}
      \mu_{i} \UE{\dot{\sigma}}_{i}(t) 
      &\mathrm{if}\quad  t \in \UE{\ClT}_{i}
      \\
      \sum_{j;j\geq i+1} 
      \sum_{k \in \ClK}
      \int_{h \in \ClH} h \UE{q}_{j,k}(h,t) \mathrm{d}h
      &\mathrm{if}\quad  t \notin \UE{\ClT}_{i}
    \end{dcases}
    &&\forall i \in \ClI, \quad \forall t \in \ClT.
  \end{align}
  Because the following inequality holds:
  \begin{align}
    \sum_{j;j\geq i+1} 
    \sum_{k \in \ClK}
    \int_{h \in \ClH} h \UE{q}_{j,k}(h,t) \mathrm{d}h
    \leq \mu_{i+1} < \mu_{i}.   
  \end{align}
  we confirm that the queueing condition \eqref{eq:DUE_QueueingCondition} is satisfied.
  \par
  Subsequently, we check the departure time choice condition \eqref{eq:DUE_DeparturePrefferd_Time_Choice}.
  Because $\UE{\ClT}_{i} = \SO{\ClT}_{i}$, $\UE{w}_{i}(t) = \SO{p}_{i}(t)$ and $\UE{\lambda}_{i}(h) = \SO{\lambda}_{i}(h)$, the following condition holds:
  \begin{align}
    &\UE{\ClT}_{i} \cap \widehat{\ClT}_{k}
    = 
    \left\{
      t \in \ClT \mid  
    h \left( c_{k}(t) + \sum_{j;j \leq i} \UE{w}_{j}(t) \right)= \UE{\lambda}_{i}(h)
    \right\}
    && 
    \forall t \in \ClT, \quad \forall k \in \ClK.
  \end{align}
  Thus, if we guarantee the following non-negativity condition, the departure time choice condition \eqref{eq:DUE_DeparturePrefferd_Time_Choice} is satisfied:
  \begin{align}
    t \in \UE{\ClT}_{i} \cap \widehat{\ClT}_{k}
    \quad 
    \Rightarrow
    \quad
    \int_{h \in \ClH} \UE{q}_{i,k}(h,t) \mathrm{d}h > 0.
  \end{align}
  Because $h > 0$, $\forall h \in \ClH$, this inequality is equivalent to the following inequality:
  \begin{align}
    F_{i,k}(t)
    \equiv \int_{h \in \ClH} h \UE{q}_{i,k}(h,t) 
    \mathrm{d}h > 0.
  \end{align}
  Consider $i \in \ClI \setminus \{ I \}$, 
  $F_{i,k}(t)$ is calculated as follows:
  \begin{align}
    F_{i,k}(t) &= 
    \left(
      \mu_{i} \UE{\dot{\sigma}}_{i}
       - 
      \mu_{i+1} \UE{\dot{\sigma}}_{i+1}
    \right)  
    \\
    &= \begin{dcases}
   \left( \mu_{i}-\mu_{i+1} \right)
    \left( 1 + \dot{c}_{k}(t)\right) 
    &\mathrm{if}\quad  t \in \UE{\ClT}_{i-1}
    \\
      \mu_{i} 
       - 
      \mu_{i+1} \left( 1 + \dot{c}_{k}(t) \right)
    &\mathrm{if}\quad  t \notin \UE{\ClT}_{i-1}
    \end{dcases}
    &&\forall t \in \UE{\ClT}_{i} \cap \widehat{\ClT}_{k}.
  \end{align}
  Because the schedule delay cost function satisfies the following condition:
  \begin{align}
   &- 1 < \dot{c}_{k}(t)
   &&\forall k \in \ClK, \quad \forall t \in \ClT,
   \\
   &\dot{c}_{k}(t) < \dfrac{\mu_{i} - \mu_{i+1}}{\mu_{i+1}}
   &&\forall k \in \ClK, \quad \forall t \in \ClT,
  \end{align}
  we find $F_{i,k}(t) > 0$.
  \par
  In the case of $i=I$, $F_{I,k}(t)$ is calculated as follows:
  \begin{align}
    F_{I,k}(t) 
    &= \mu_{I} \UE{\dot{\sigma}}_{I}
    = \mu_{I} \left( 1 + \dot{c}_{k}(t) \right)
    &&\forall t \in \UE{\ClT}_{I} \cap \widehat{\ClT}_{k}.
  \end{align}
  Because $- 1 < \dot{c}_{k}(t)$, $F_{I,k}(t)>0$, i.e., the non-negativity condition is satisfied.
  \par
  Consequently, the variables shown in \Cref{pro:DUE_Solution} satisfy all short-term equilibrium conditions under the QRP condition.
  \hfill \Halmos

\subsection{Proof of \Cref{lem:lambda_NSvsSWH}}
\label{subsec:prf_lem_lambda_NSvsSWH}
We calculate $\NS{\lambda}_{i}(1) - \SWH{\lambda}_{i}(1)$ as follows:
\begin{align}
  &\NS{\lambda}_{i}(1) - \SWH{\lambda}_{i}(1) 
  = 
    \overline{c}_{1}(A_{i}, \overline{\mu}_{i})
   -
    \overline{c}_{2}(A_{i}, \overline{\mu}_{i})
  = d \delta > 0
    &&\forall i \in \ClI.
\end{align}
Thus, we obtain \eqref{eq:lambda_NSvsSWH}. 
\hfill \Halmos

\subsection{Proof of \Cref{lem:r_NSvsSWH}}
\label{subsec:prf_lem_r_NSvsSWH}
We calculate $\Delta r_{i} \equiv r_{i} - r_{i+1}$ in scenarios (NS) and (SWH) as follows:
\begin{align}
  \Delta \NS{r}_{i}
  &= 
    \left\{ \theta^{O} - 
    \overline{c}_{1}(A_{i}, \overline{\mu}_{i}) - \sum_{j;j\leq i} f_{j} + \NS{\rho} 
    \right\}
   - 
    \left\{ \theta^{O} - 
    \overline{c}_{1}(A_{i+1}, \overline{\mu}_{i+1}) - \sum_{j;j\leq i+1} f_{j} + \NS{\rho}
    \right\}
   \notag
   \\
   & = - \dfrac{A_{i}}{\overline{\mu}_{i}}\delta 
       + \dfrac{A_{i+1}}{\overline{\mu}_{i+1}}\delta 
       + f_{i+1}
       \qquad \qquad \qquad \forall i \in \ClI \setminus \{ I \},
  \\
  \Delta \SWH{r}_{i} &= 
    \left\{ \theta^{O} - 
    \overline{c}_{2}(A_{i}, \overline{\mu}_{i}) - \sum_{j;j\leq i} f_{j} + \SWH{\rho} 
    \right\}
     -
    \left\{ \theta^{O} - 
    \overline{c}_{2}(A_{i+1}, \overline{\mu}_{i+1}) - \sum_{j;j\leq i+1} f_{j} + \SWH{\rho}
    \right\}
   \notag
   \\
   & = - \dfrac{A_{i}}{\overline{\mu}_{i}}\delta 
       + d \delta 
       + \dfrac{A_{i+1}}{\overline{\mu}_{i+1}}\delta 
       - d \delta
       + f_{i+1} 
   = 
       - \dfrac{A_{i}}{\overline{\mu}_{i}}\delta 
       + \dfrac{A_{i+1}}{\overline{\mu}_{i+1}}\delta 
       + f_{i+1}
       \quad \quad  \forall i \in \ClI \setminus \{ I \}.
\end{align}
Hence, we have $\Delta \NS{r}_{i} = \Delta \SWH{r}_{i}$.
Based on this, we obtain \eqref{eq:r_NSvsSWH} under boundary conditions $\NS{r}_{I} = 0$ and $\SWH{r}_{I} = 0$.
This completes the proof.
\hfill \Halmos

\subsection{Proof of \Cref{lem:lambda_NSvsTLC}}
\label{subsec:prf_lem:lambda_NSvsTLC}
The workers who live at the location $i<i^{\ast}$ choose the office work ratio $h=1$ even if TLC is introduced.
Thus, \eqref{eq:lambda_NSvsTLC_1} holds because their commuting cost does not change before and after the introduction of TLC.
\par
The equilibrium commuting at the location $i > i^{\ast}$costs decrease introducing remote work zone at the location $i > i^{\ast}$ because workers lived in there do not commute to the CBD.
The workers who live in the mixed zone adjust their office work ratio so that their utility is not lower than one from TLC.
Thus, \eqref{eq:lambda_NSvsTLC_2} holds.
\hfill \Halmos

\subsection{Proof of \Cref{lem:r_NSvsTLC}}
\label{subsec:prf_lem:r_NSvsTLC}
In the remote work zone and the mixed zone, the equilibrium land rent equal to one in the most upstream location, i.e., $\TLC{r}_{i}=\TLC{r}_{I}=0$, $i \geq \TLC{i}{}^{\ast}$.
\par
In the office work zone, $\Delta r_{i} \equiv r_{i} - r_{i+1}$ in scenarios (NS) and (Rx) as follows:
\begin{align}
  \Delta \NS{r}_{i} &= 
    \left\{ \theta^{O} - 
    \overline{c}_{1}(A_{i}, \overline{\mu}_{i}) - \sum_{j;j\leq i} f_{j} + \NS{\rho} 
    \right\}
    - 
    \left\{ \theta^{O} - 
    \overline{c}_{1}(A_{i+1}, \overline{\mu}_{i+1}) - \sum_{j;j\leq i+1} f_{j} + \NS{\rho}
    \right\}
   \notag
   \\
   & = - \dfrac{A_{i}}{\overline{\mu}_{i}}\delta 
      + \dfrac{A_{i+1}}{\overline{\mu}_{i+1}}\delta 
       + f_{i+1}
       \qquad \qquad \qquad  \forall i < i^{\ast},
  \\
  \Delta \TLC{r}_{i}
   &= 
    \left\{ \theta^{O} - 
    \overline{c}_{1}(A_{i}, \overline{\mu}_{i}) - \sum_{j;j\leq i} f_{j} + \SWH{\rho} 
    \right\}
    - 
    \left\{ \theta^{O} - 
    \overline{c}_{1}(A_{i+1}, \overline{\mu}_{i+1}) - \sum_{j;j\leq i+1} f_{j} + \SWH{\rho}
    \right\}
   \notag
   \\
   & =  - \dfrac{A_{i}}{\overline{\mu}_{i}}\delta 
   + \dfrac{A_{i+1}}{\overline{\mu}_{i+1}}\delta 
       + f_{i+1}
   \qquad \qquad \qquad \forall i < i^{\ast}.
\end{align}
Hence, we have $\Delta \NS{r}_{i} = \Delta \TLC{r}_{i}$, $i<i^{\ast}$.
Based on this, we obtain \eqref{eq:r_NSvsTLC} because $\NS{r}_{i} > \TLC{r}_{\TLC{i}{}^{\ast}} = 0$.
This completes the proof.

\subsection{Proof of \Cref{lem:lambda_SWHvsCS}}
\label{subsec:prf_lem:lambda_SWHvsCS}
\Cref{lem:lambda_SWHvsCS} is proved by the same way as \Cref{lem:lambda_NSvsTLC}.
\hfill \Halmos

\subsection{Proof of \Cref{lem:r_SWHvsCS}}
\label{subsec:prf_lem:r_SWHvsCS}
\Cref{lem:r_SWHvsCS} is proved by the same way as \Cref{lem:r_NSvsTLC}.
\hfill \Halmos

\subsection{Proof of \Cref{lem:lambda_TLCvsCS}}
\label{subsec:prf_lem:lambda_TLCvsCS}
As shown in \Cref{tab:summary}, the changes in the type of location when introducing SWH in addition to TLC are classified into five cases.
\par
In the first case (i.e., office work zone $\rightarrow$ office work zone), the commuting cost decreases because the following inequality holds:
\begin{align}
  \overline{c}_{1}(X, \mu) > \overline{c}_{2}(X, \mu).
\end{align}
In the second case (i.e., mixed zone $\rightarrow$ office work zone), the commuting cost also decreases.
In the third case (i.e., remote work zone $\rightarrow$ office work zone) and forth case (i.e., remote work zone $\rightarrow$ mixed zone), the commuting cost increases because the workers who live in there start to the commute to the CBD with introducing SWH.
In the fifth case (i.e., mixed zone $\rightarrow$ mixed zone), the commuting cost does not change.
Thus, we obtain \Cref{lem:lambda_TLCvsCS}.
\hfill \Halmos

\subsection{Proof of \Cref{lem:r_TLCvsCS}}
\label{subsec:prf_lem:r_TLCvsCS}
We calculate $\Delta \CS{r}_{i} \equiv \CS{r}_{i} - \CS{r}_{i+1}$ in scenarios (CS) as follows:
\begin{align}
  \Delta \CS{r}_{i} &= 
   \left\{ \theta^{O} - 
   \overline{c}_{2}(A_{i}, \overline{\mu}_{i}) - \sum_{j;j\leq i} f_{j} + \SWH{\rho} 
   \right\}
    - 
   \left\{ \theta^{O} - 
   \overline{c}_{2}(A_{i+1}, \overline{\mu}_{i+1}) - \sum_{j;j\leq i+1} f_{j} + \SWH{\rho}
   \right\}
  \notag
  \\
  & =  - \dfrac{A_{i}}{\overline{\mu}_{i}}\delta 
  + \dfrac{A_{i+1}}{\overline{\mu}_{i+1}}\delta 
      + f_{i+1}
      \qquad \qquad \forall i < \CS{i}{}^{\ast}.
\end{align}
By comparing $\Delta \CS{r}_{i}$ with $\Delta \TLC{r}_{i}$, we obtain \Cref{lem:r_TLCvsCS} because $\CS{i}{}^{\ast} > \TLC{i}{}^{\ast}$ and $\CS{r}_{\CS{i}{}^{\ast}}=\TLC{r}_{\TLC{i}{}^{\ast}} = 0$.

\subsection{Proof of \Cref{thm:SWH_Effect,thm:RemoteWorkEffect,thm:SWHvsCS,thm:RemoteWorkParadox}}
\label{subsec:prf_thm}
This section proves \Cref{thm:SWH_Effect,thm:RemoteWorkEffect,thm:SWHvsCS,thm:RemoteWorkParadox} simultaneously.
\par
We first derive the inequalities related to equilibrium utility as follows:
\begin{align}
  &\SWH{\rho} - \NS{\rho} 
  =
  \theta^{O} - \dfrac{R_{I}}{\mu} \delta + d \delta
   - \sum_{i \in \ClI} f_{i}
   - \theta^{O} + \dfrac{R_{I}}{\mu} \delta 
   +  \sum_{i \in \ClI} f_{i} 
   = d \delta > 0,
  \\
  &\TLC{\rho} - \NS{\rho}
 = \theta^{R}
 - \theta^{O} + \dfrac{R_{I}}{\overline{\mu}_{I}} \delta + \sum_{i \in \ClI} f_{i}
 > 0,
 \\
 &\CS{\rho} - \SWH{\rho}
 = \theta^{R}
 - \theta^{O} + \dfrac{R_{I}}{\overline{\mu}_{I}} \delta - d \delta + \sum_{i \in \ClI} f_{i}
 > 0,
 \\
 &\CS{\rho} - \TLC{\rho}
  = \theta^{R} - \theta^{R} = 0.
\end{align}
Subsequently, we confirm the inequalities related to total commuting cost.
The total commuting cost in each scenario can be calculated as follows:
\begin{align}
  &\NS{TC} =   \sum_{i \in \ClI} \int_{t \in \Gamma_{1}(\overline{c}_{1}(A_{i}, \overline{\mu}_{i}))} 
  \overline{c}_{1}(A_{i}, \overline{\mu}_{i}) \overline{\mu}_{i} \mathrm{d} t,
  \\ 
  &\SWH{TC} =    \sum_{i \in \ClI} \int_{t \in \Gamma_{2}(\overline{c}_{2}(A_{i}, \overline{\mu}_{i}))} 
  \overline{c}_{2}(A_{i}, \overline{\mu}_{i}) \overline{\mu}_{i} \mathrm{d} t,
  \\
  &\TLC{TC} = \sum_{i \in \ClI} \int_{t \in \Gamma_{1}(\overline{c}_{1}(\TLC{\eta}(i)A_{i}, \overline{\mu}_{i}))} 
  \overline{c}_{1}(\TLC{\eta}(i)A_{i}, \overline{\mu}_{i}) \overline{\mu}_{i} \mathrm{d} t,
  \\
  &\CS{TC}  = \sum_{i \in \ClI} \int_{t \in \Gamma_{2}(\overline{c}_{2}(\CS{\eta}(i)A_{i}, \overline{\mu}_{i}))} 
  \overline{c}_{2}(\CS{\eta}(i)A_{i}, \overline{\mu}_{i}) \overline{\mu}_{i} \mathrm{d} t.
\end{align}
By definition of the function $\Gamma_{K}(c)$, we have following properties:
\begin{align}
  &\| \Gamma_{K}(c) \| > \| \Gamma_{K+1}(c) \|
  &&\forall c \in (0, +\infty),
  \\
  &\| \Gamma_{K}(c) \| < \| \Gamma_{K}(c+\Delta c) \|
  &&\forall K = 1,2.
\end{align}
In addition, the function $\overline{c}_{K}(X, \mu)$ satisfies the following inequality:
\begin{align}
  \overline{c}_{1}(X, \mu) > \overline{c}_{2}(X, \mu).
\end{align}
  From these properties, we obtain $\SWH{TC} < \NS{TC}$, $\TLC{TC} < \NS{TC}$, and $\CS{TC} < \SWH{TC}$.
  Moreover, \Cref{sec:ExampleoftheParadox} also shows an example of a case where $\CS{TC} > \TLC{TC}$.
  This completes the proof.~\hfill \Halmos

\section{Example of the Paradox}\label{sec:ExampleoftheParadox}
  This section shows an example of anlytical solution in each scenario.
  We numerically confirm the analysis result in Section 5 and show that the paradox shown in \Cref{thm:RemoteWorkParadox} can occur.
  We consider the cases of $N=3$ and $K=2$.
  The parameters are chosen as follows:
\begin{align}
  &\Vtmu = [70, 40, 10], \quad 
  \VtA = [750, 1500, 700], 
  \quad \Vtf = [1.5, 1.0, 1.0], 
  \quad \ClT = [0, 100],
  \notag
  \\
  &\beta = 0.30, \quad 
  \gamma = 0.60, \quad 
  \theta^{O} = 40.0, \quad 
  \theta^{R} = 30.0,
  \notag
  \\
  &t_{1} = 60 \text{ in the scenarios (NS) and (TLC)},
  \quad t_{1}, t_{2} = 50, 70 \text{ in the scenarios (SWH) and (CS)}.
  \notag
\end{align}
The parameters satisfy the QRP condition \eqref{eq:QRP_Condition}.
\par
  First, we show the location patterns and the equilibrium land rent pattern in \Cref{tab:example_location_pattern}.
    The table shows that the location $i=2$ is the mixed zoned in scenario (TLC), but is the office work zone in scenario (CS).
  The office work ratio in the mixed zone is calculated by using \eqref{eq:eta_i}.
\par
  Subsequently, we show the equilibrium commuting costs, total commuting costs and equilibrium utilities in \Cref{tab:example_equilibrium}.
  We observe that equilibrium utility in the scenario (CS) is equal to that in the scenario (TLC); however, the total commuting cost in the scenario (CS) is higher than that in the scenario (TLC).
  This is the statement of \Cref{thm:RemoteWorkParadox}.
\begin{table}[t]
  \caption{Equilibrium Location Pattern in each Scenario.}
  \label{tab:example_location_pattern}
  \centering
  \renewcommand{\arraystretch}{1.5} %
  \begin{tabular}{|c||c|c|c|}
    \hline
     {} & $i=3$ & $i=2$ & $i=1$ \\
     \hline \hline
    NS & 
    \begin{tabular}{c}
      Office work zone \\
      $\NS{r}_{3} = 0.0$
    \end{tabular}
     & 
     \begin{tabular}{c}
      Office work zone \\ 
      $\NS{r}_{2} = 5.0$ 
      \end{tabular} 
     & 
     \begin{tabular}{c} 
      Office work zone \\ 
      $\NS{r}_{1} = 11.0$
      \end{tabular} \\
    \hline
    SWH & \begin{tabular}{c}
      Office work zone \\
      $\SWH{r}_{3} = 0.0$
    \end{tabular}
     & 
     \begin{tabular}{c}
      Office work zone \\ 
      $\SWH{r}_{2} = 5.0$
      \end{tabular} 
     & 
     \begin{tabular}{c} 
      Office work zone \\ 
      $\SWH{r}_{1} = 11.0$
      \end{tabular} \\
    \hline
    TLC & 
    \begin{tabular}{c}
      Remote work zone \\
      $\TLC{r}_{3} = 0.0$
    \end{tabular}
      & 
      \begin{tabular}{c}
        Mixed zone, $\TLC{\eta}(2)=0.75$\\ 
        $\TLC{r}_{2} = 0.0$
        \end{tabular} 
      & 
      \begin{tabular}{c} 
        Office work zone \\ 
        $\TLC{r}_{1} = 3.5$
        \end{tabular} \\
    \hline
    CS &
    \begin{tabular}{c}
      Mixed zone, $\CS{\eta}(3)=0.75$\\
      $\CS{r}_{3} = 0.0$
    \end{tabular}
    &
    \begin{tabular}{c}
      Office work zone\\ 
      $\CS{r}_{2} = 1.5$
      \end{tabular}
    &
    \begin{tabular}{c} 
      Office work zone \\ 
      $\CS{r}_{1} = 7.5$
      \end{tabular} \\
    \hline
  \end{tabular}
\end{table}
\begin{table}[t]
  \caption{Equilibrium Commuting Cost, Total Commuting Cost and Equilibrium Utility in each Scenario.}
  \label{tab:example_equilibrium}
  \centering
  \renewcommand{\arraystretch}{1.5} %
  \begin{tabular}{|c||c|c|c||r||r|}
    \hline
     {} & $i=3$ & $i=2$ & $i=1$ 
     & Total commuting cost & 
       Equilibrium utility \\
    \hline
    \hline
    NS & $\NS{\lambda}_{3} = 14.0$ & $\NS{\lambda}_{2} = 10.0$ & $\NS{\lambda}_{1} = 5.0$ 
    & $\NS{TC}= 28550.0$ & $\NS{\rho}= 22.5$
    \\
    \hline
    SWH & $\SWH{\lambda}_{3} = 10.0$ & $\SWH{\lambda}_{2} = 6.0$ & $\SWH{\lambda}_{1} = 1.0$ 
    & $\SWH{TC}= 16750.0$ & $\SWH{\rho}= 26.5$ \\
    \hline
    TLC & $\TLC{\lambda}_{3} = 0.0$ & $\TLC{\lambda}_{2} = 7.5 $ & $\TLC{\lambda}_{1} = 5.0$ 
    & $\TLC{TC}= 12187.5$ & $\TLC{\rho}= 30.0$ \\
    \hline
    CS &$\CS{\lambda}_{3} = 6.6$ & $\CS{\lambda}_{2} = 6.0$ & $\CS{\lambda}_{1} = 1.0$ 
    & $\CS{TC}= 13162.5$ & $\CS{\rho}= 30.0$ \\
    \hline
  \end{tabular}
\end{table}

\bibliographystyle{apalike} 
\bibliography{reference}
\end{document}